\documentclass[12pt, english]{amsart}

\usepackage{stmaryrd}

\usepackage{amscd,amsxtra}
\usepackage[all]{xy}
\usepackage{babel}
\usepackage{amstext}
\usepackage{amsmath}
\usepackage{amsfonts}

\usepackage{xypic}
\xyoption{all}

\usepackage[dvips]{graphics}
\usepackage{epsfig,graphicx}
\usepackage{amscd}

\usepackage{color, colortbl}
\usepackage{psfrag}

\input xy

\DeclareMathOperator{\rk}{rk}
\DeclareMathOperator{\Coh}{\mathsf{Coh}}
\DeclareMathOperator{\VB}{\mathsf{VB}}
\DeclareMathOperator{\TF}{\mathsf{TF}}
\DeclareMathOperator{\BM}{\mathsf{BM}}

\DeclareMathOperator{\Pic}{Pic}

\DeclareMathOperator{\Hom}{\mathsf{Hom}}
\DeclareMathOperator{\Ext}{\mathsf{Ext}}

\DeclareMathOperator{\Aut}{\mathsf{Aut}}
\DeclareMathOperator{\End}{\mathsf{End}}
\DeclareMathOperator{\Spec}{\mathsf{Spec}}

\newcommand{\lar}{\longrightarrow}
\newcommand{\kF}{\mathcal F}

\newcommand{\kG}{\mathcal G}
\newcommand{\kL}{\mathcal L}
\newcommand{\kB}{\mathcal B}
\newcommand{\kM}{\mathcal M}
\newcommand{\kO}{\mathcal O}
\newcommand{\kP}{\mathcal P}
\newcommand{\kI}{\mathcal I}

\newcommand{\kJ}{\mathcal J}
\newcommand{\kE}{\mathcal E}
\newcommand{\kT}{\mathcal T}
\newcommand{\kN}{\mathcal N}
\newcommand{\kS}{\mathcal S}
\newcommand{\kk}{{\bf k}}
\newcommand{\EE}{\mathbb{E}}
\newcommand{\TT}{\mathbb{T}}
\newcommand{\XX}{\mathbb{X}}
\newcommand{\PP}{\mathbb{P}}
\newcommand{\ZZ}{\mathbb{S}}
\newcommand{\LL}{\mathbb{L}}
\newcommand{\YY}{\mathbb{Y}}
\renewcommand{\SS}{\mathbb{S}}

\newcommand{\kA}{\mathcal A}
\newcommand{\mZ}{\mathbb Z}
\newcommand{\mP}{\mathbb P}
\newcommand{\mA}{\mathbb A}

\renewcommand{\mod}{\mathop{\rm mod}}

\newcommand{\vr}{\vrule width 0.1mm }
\newcommand{\fvr}{\vrule width 0.3mm }
\newcommand{\mk}[1]{\makebox[15pt][c]{\rule{0pt}{12pt}~$ #1 $ }}
\newcommand{\bmk}[1]{\makebox[15pt][c]{\rule{0pt}{12pt}~$ #1 $ }}
\newcommand{\dhr}{
\begin{picture}(10,0)
\thinlines
\multiput(-5,3)(5,0){5}{\circle*{1}}
\end{picture}
}

\newtheorem{theorem}{Theorem}[section]
\newtheorem{algorithm}[theorem]{Algorithm}
\newtheorem{proposition}[theorem]{Proposition}
\newtheorem{lemma}[theorem]{Lemma}
\newtheorem{corollary}[theorem]{Corollary}

\theoremstyle{definition}
\newtheorem{definition}[theorem]{Definition}
\newtheorem{problem}[theorem]{Problem}
\newtheorem{remark}[theorem]{Remark}
\newtheorem{example}[theorem]{Example}
\begin{document}

\title
[Vector bundles on degenerations of elliptic curves]
{Vector bundles and torsion free sheaves on degenerations of elliptic curves}

\author{Lesya Bodnarchuk}
\address{%
Technische Universit\"at Kaiserslautern,
Erwin-Schr\"odinger-Stra\ss{}e\\
67663 Kaisers\-lautern,
 Germany}
\email{bodnarch@mathematik.uni-kl.de}

\author{Igor Burban }
\address{%
Institut f\"ur Mathematik,
Johannes Gutenberg-Universit\"at Mainz,
55099 Mainz, Germany
}
\email{burban@mathematik.uni-mainz.de}

\author{Yuriy Drozd} 
\address{%
 Kyiv Taras Shevchenko University,
Department of Mechanics and Mathematics
Volo\-dimirska 64,
01033 Kyiv, Ukraine
}
\email{yuriy@drozd.org}

\author{Gert-Martin Greuel}
\address{%
Technische Universit\"at Kaiserslautern,
Erwin-Schr\"odinger-Stra\ss{}e,
67663 Kaisers\-lautern,
 Germany
 }
\email{greuel@mathematik.uni-kl.de}

\maketitle

%%%%%%%%%%%%%%%%%
%Your text goes here. Separate text sections with the standard \LaTeX\ sectioning commands.

\begin{abstract}
In this paper we give a survey about the classification of  vector bundles
and torsion free  sheaves on degenerations of elliptic curves.
Coherent sheaves on singular curves of arithmetic genus one can be studied using the technique
of matrix problems or via Fourier-Mukai transforms, both methods are discussed here.
 Moreover, we include new proofs of some classical results about vector bundles on elliptic curves.
\end{abstract}

\tableofcontents

\section{Overview}

The aim of this paper  is to give a survey  on  the classification of vector
bundles and
torsion free  sheaves on
singular projective curves of arithmetic genus one.  We include new proofs of
some classical results on coherent sheaves on smooth elliptic curves, which use  the technique
of derived categories and Fourier-Mukai  transforms and are simpler than the original ones.
Some results about singular curves are new or at least presented in a new framework.

This research  project
had several sources of   motivation and inspiration.
Our study of vector bundles on degenerations of elliptic curves was originally motivated by
the
McKay correspondence  for minimally elliptic surface singularities
\cite{Kahn}.    Here we use as the main technical tool  methods
from the representation theory of associative algebras, in particular, a  key tool in our approach to classification problems
 is played by  the technique of
``representations of bunches of chains'' or ``Gelfand problems'' \cite{mp}.
At last, but
not least we want   to mention that our  research was strongly influenced by ideas and methods
coming from  the
homological mirror symmetry \cite{Kontsevich, PolishchukZaslow, FMW}.

 Many different questions concerning  properties of the category of
vector bundles and coherent sheaves on degenerations of elliptic curves are  encoded
in the  following general set-up:

\begin{figure}[ht]
\hspace{1.5cm}
\psfrag{a}{0}
\psfrag{b}{$t$}
\psfrag{A}{$D^b(\Coh(\EE_0))$}
\psfrag{B}{$D^b(\Coh(\EE_t))$}
\includegraphics[height=4.5cm,width=6.0cm,angle=0]{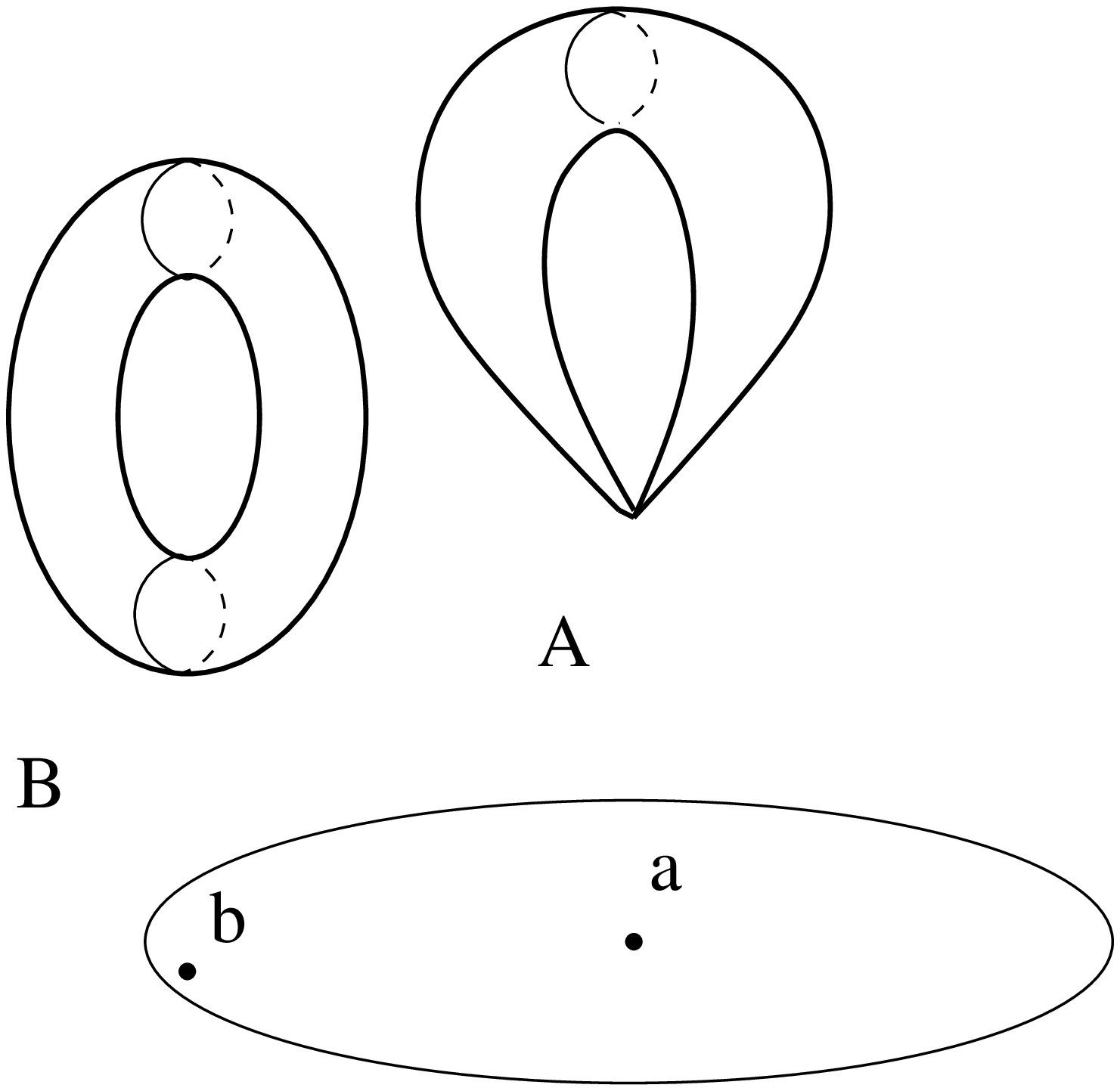}
\end{figure}

\begin{problem}
Let $\EE \lar \TT$ be a flat family of
projective curves of arithmetic genus one such that the
fiber $\EE_t$ is smooth for generic $t$   and singular for $t = 0$.
What happens with the derived category
$D^b(\Coh_{\EE_t})$, when $t\rightarrow  0$?
\end{problem}

In order to start working on this  question one has to consider  the absolute case first,
where
the base $\TT$ is $\Spec(\kk)$. In particular, one has to describe
indecomposable
vector bundles and indecomposable objects of the derived category of
coherent sheaves on degenerations of elliptic curves
and develop a technique to
calculate
homomorphism and extension spaces between indecomposable torsion free sheaves
 as well as various operations on them,    like tensor products and  dualizing.

For the first time we  face
this sort of    problems  when dealing
 with  the  McKay correspondence  for \emph{minimally elliptic singularities}.
Namely, let $\mathbb{S} = \Spec({\bf R })$
be the spectrum  of a complete (or analytical)
two-dimensional minimally elliptic singularity,
$\pi: \widetilde{\XX} \lar \mathbb{S}$ its   minimal resolution,  and
$\EE$ the
exceptional divisor.   Due to a construction of Kahn \cite{Kahn},
the functor $M \mapsto
res_{\EE}(\pi^*(M)^{\vee\vee})$ establishes
a bijection   between   the reflexive
${\bf R}$--modules  (maximal Cohen-Macaulay modules) and the
generically globally generated indecomposable vector bundles on $\EE$ with vanishing
first cohomology.\footnote{$res_{\EE}$ denotes the restriction to $\EE$}

A typical example of a minimally elliptic singularity is a
$T_{pqr}$\,--\,singularity,  given  by the equation
$\kk\llbracket x,y,z\rrbracket  /(x^p + y^q + z^r - \lambda xyz)$, where
$\frac{1}{p} + \frac{1}{q}+  \frac{1}{r} \le   1$ and $\lambda \ne 0$.
If $\frac{1}{p} + \frac{1}{q} + \frac{1}{r} = 1$, then this singularity is
\emph{simple elliptic} and the exceptional divisor $\EE$ is a smooth elliptic curve.
Thus, in this case a   description  of indecomposable maximal Cohen-Macaulay modules follows
from  Atiyah's classification of
vector bundles on elliptic curves \cite{Atiyah}.
The main result of Atiyah's  paper  essentially says:

\begin{theorem}[Atiyah]
An indecomposable vector bundle $\kE$ on an elliptic curve $\EE$ is uniquely determined by its
rank $r$, degree $d$  and  determinant $det(\kE) \in \Pic^d(\EE)\cong \EE$.
\end{theorem}

In Section \ref{subsec:ellcurves} we give a new proof of this result.
However, if $\frac{1}{p} + \frac{1}{q}+  \frac{1}{r} <  1$,
then $\mathbb{S}$ is a  so-called \emph{cusp   singularity} and in this case
$\EE$ is  a cycle
 of $n$  projective lines $\EE_n$, where $\EE_1$ denotes  a rational curve with one node.

%%%%%%%
%my figure:
\begin{figure}[ht]
%\hspace{0.7cm}
\begin{minipage}[b]{3cm}
\includegraphics[height=2.5cm,width=2cm,angle=0]{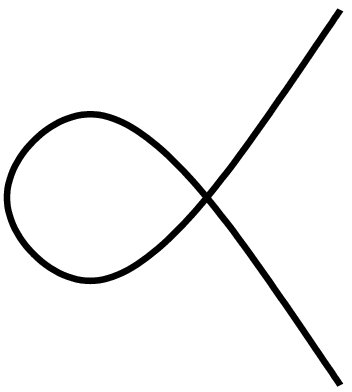}
\center {$\EE_1\phantom{BB}$}
%\caption*{c)}
\end{minipage}
\hspace{1cm}
\begin{minipage}[b]{3cm}
\hspace{0.3cm}
\includegraphics[height=2.5cm,width=2.5cm,angle=0]{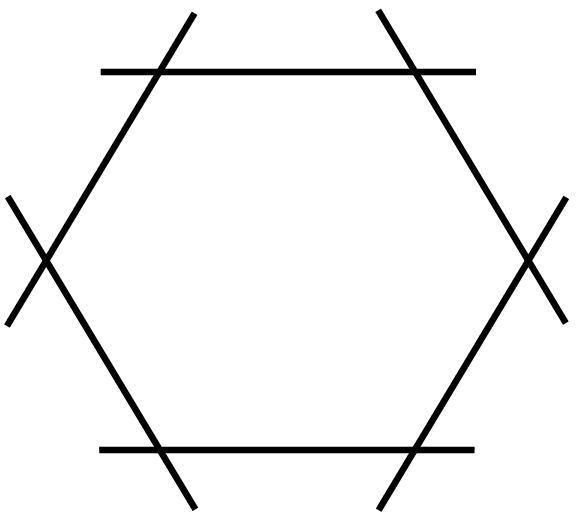}
\center {$\phantom{Bb}\EE_6$}
%\caption*{c)}
\end{minipage}
\end{figure}

A complete classification  of indecomposable vector bundles and torsion free sheaves
on these curves in
the case of an   arbitrary base field $\kk$
was  obtained by  Drozd and Greuel
\cite{DrozdGreuel}. For algebraically closed fields
there is the following description, which we prove  in  Section \ref{subsec:prop}.

\begin{theorem}\label{thm:classification}
Let $\EE_n$ denote a cycle of $n$ projective lines and $\mathbb{I}_k$  be a  chain of
$k$ projective lines, $\kE$ an indecomposable
torsion free sheaf on $\EE_n$.
\begin{enumerate}
\item If $\kE$ is  locally free, then there is an \'etale covering
$\pi_r: \EE_{nr}\lar \EE_n$, a line bundle $\kL\in \Pic(\EE_{nr})$ and
a natural number $m\in {\mathbb N}$ such that
$$
 \kE \cong \pi_{r*}(\kL\otimes \kF_m),
$$
where $\kF_m$ is an indecomposable vector bundle on $\EE_{nr}$, recursively defined
by the sequences
$$
0 \lar \kF_{m-1} \lar \kF_m \lar \kO \lar 0, \quad  m\ge 2, \quad \kF_1 = \kO.
$$
\item If $\kE$ is not locally free then there exists a finite  map
$p_k: \mathbb{I}_{k} \lar \EE_n$ and a line bundle $\kL \in \Pic(\mathbb{I}_k)$ (where
$k$, $p_k$ and  $\kL$ are determined by $\kE$) such that
$\kE \cong p_{k*}(\kL)$.
\end{enumerate}
\end{theorem}

This classification is completely analogous  to   Oda's description of vector bundles
on smooth elliptic curves \cite{Oda} and   provides quite simple rules for the computation
of the decomposition of the tensor product of any two vector bundles into
a direct sum of indecomposable ones. It allows to describe  the dual sheaf  of an
indecomposable torsion free sheaf  as well as the  dimensions of homomorphism and extension
spaces between indecomposable vector bundles  (and in particular, their cohomology), see
\cite{Burban, BDG}. We carry out these computations in Section \ref{subsec:prop}.

However, the way we prove this theorem,  essentially uses ideas from
representation theory and the technique of matrix problems \cite{mp}.
Using a similar approach,
Theorem \ref{thm:classification}  was generalized by Burban and Drozd \cite{BurbanDrozd1}
to classify
indecomposable
complexes of the bounded (from the right) derived category of coherent
sheaves $D^-(\Coh(\EE))$ on
a cycle of projective lines $\EE = \EE_n,$
see also \cite{BurbanDrozd2} for the case of associative algebras.

The situation turns out to be quite different for other
singular projective curves
of arithmetic genus one.
For example, in the case of a cuspidal rational curve $zy^2 = x^3$
even a classification of indecomposable semi-stable
vector bundles of a given slope  is a \emph{representation-wild}
problem.\footnote{
An exact $\kk$--linear category $\mathsf{A}$ over an algebraically closed field $\kk$
 is called wild if it contains as a full
subcategory the category of finite-dimensional representations of
\emph{any} associative algebra.}
 However, if we  restrict
our attention only on  stable vector bundles,
then this problem  turns out to be \emph{tame}
again.\footnote{For a  formal definition of tameness we refer to \cite{DrozdGreuel},
where a wild-tame
dichotomy for vector bundles and torsion free sheaves on reduced curves was proven.}
 Moreover, the combinatorics of the answer is essentially the same as for
smooth and nodal Weierstra\ss{} curves\footnote{In this paper  we call a plane
cubic curve
{\it Weierstra\ss{} curve}. If $\kk$ is algebraically closed
and $char(\kk)\ne 2,3$ then it can be written in the form
$zy^2 = 4x^3 - g_2xz^2 - g_3z^3$, where $g_2, g_3 \in \kk$.
It is singular if and only if $g_2^3 = 27 g_3^2$ and unless $g_2 = g_3 =0$ the singularity
is a node, whereas in the case $g_2 = g_3 =0$ the singularity is a cusp.}:

%my change:
%\begin{theorem}[see \cite{BodnarchukDrozd, BurbanKreussler3}]
%Let $\EE$ be a cuspidal rational
% curve and $\kE$ a stable vector bundle on $\EE$. Then
%$\kE$ is completely determined by its rank $rk(\kE)$, its degree $deg(\kE)$
% and its determinant $det(\kE) \in \Pic^{\deg(\kE)}(\EE) \cong \kk$.
%Moreover,  $rk(\kE)$ and $deg(\kE)$ are mutually prime.
%\end{theorem}

 \begin{theorem} [see \cite{BodnarchukDrozd, BurbanKreussler3}]
  Let $\EE$ be a cuspidal cubic curve over an algebraically closed field $\kk$ then
%\begin{itemize}
%\item
%  The rank $r$ and the degree $d$ of a stable vector bundle
%$\kE$ over $\EE$ are coprime.
%\item
a stable vector bundle $\kE$
  is completely determined  by its rank $r$, its degree $d,$ that should be coprime,
  and its determinant $det(\kE) \in \Pic^d(\EE) \cong \kk$.
% $\vb^s(r,d)=\mA^1$
%\end{itemize}
 \end{theorem}

The technique of matrix problems is a very convenient  tool for the
study of vector bundles on a given
singular projective rational curve of arithmetic genus one. However,
to investigate  the
behavior of the category of  coherent sheaves  on genus one curves
in families one needs other methods.
One possible approach is
provided by the technique of derived categories and Fourier-Mukai transforms
\cite{Mukai, SeidelThomas}, see Section \ref{subsec:FMT}.
The key idea of this method is that  we can map  a sky-scraper
sheaf into a  torsion free sheaf by applying an auto-equivalence
of the derived category.
In a relative setting of elliptic fibrations with a section  one can use \emph{relative}
 Fourier-Mukai  transforms to construct examples of relatively semi-stable torsion free sheaves, see
for example \cite{BurbanKreussler2}.

\begin{theorem}[see \cite{BurbanKreussler1}]
Let $\EE$ be an \emph{irreducible}  projective curve of arithmetic genus one over
an algebraically closed field $\kk$. Then
\begin{enumerate}
\item The group of exact auto-equivalences of the derived category
$D^b(\Coh(\EE))$ transforms  stable sheaves  into stable ones and semi-stable sheaves
into semi-stable ones.

\item For any rational number  $\nu$ the abelian category $\Coh^\nu(\EE)$
of semi-stable coherent sheaves
of slope $\nu$
is equivalent to the category $\Coh^\infty(\EE)$ of coherent torsion sheaves and
this equivalence is induced by an exact auto-equivalence of $D^b(\Coh(\EE))$.

\item For any coherent sheaf $\kF$ on $\EE$
 such that $\End(\kF) = \kk$ there exists
a point $x \in \EE$ and $\Phi \in \Aut(D^b(\Coh(\EE)))$  such that $\kF \cong
\Phi(\kk(x))$.
\end{enumerate}
\end{theorem}

This theorem shows a  fundamental  difference between a nodal and a cuspidal Weierstra\ss{}
 curve.  Namely, let $\EE$ be a singular Weierstra\ss{} curve and $s$ its singular point.
Then  the category  of finite-dimensional modules over the  complete local ring
$\widehat{\kO_{\EE,s}}$
has  different representation types  in the nodal and cuspidal cases.
For a  nodal curve,  the category of finite dimensional
representations of  $\kk\llbracket x,y\rrbracket  /xy$ is tame due to a result of Gelfand and Ponomarev
\cite{GP}. In the second case,   the category
of finite length  modules over the ring $\kk\llbracket x,y\rrbracket  /(y^2 - x^3)$ is representation wild, see
for example \cite{Drozd}.

The correspondence between sky-scraper sheaves and semi-stable vector bundles on irreducible
Weierstra\ss{}  curves was first discovered by Friedman, Morgan and Witten \cite{FMW} (see also
\cite{Teodorescu}) and afterwards  widely  used in the physical literature
under the name ``spectral cover construction''.

\begin{theorem}[see \cite{FMW}]
Let $\EE$ be an irreducible Weierstra\ss{} curve, $p_0 \in \EE$ a smooth point and $\kE$ a semi-stable
torsion free sheaf  of degree zero. Then the  sequence
$$
0 \lar H^0(\kE(p_0))\otimes \kO \stackrel{ev}\lar \kE(p_0) \lar coker(ev) \lar 0
$$
is exact. Moreover, the functor $\Phi:
\kE \mapsto coker(ev)$ establish an equivalence between
the category $\Coh^0(\EE)$ of semi-stable torsion-free sheaves of degree zero
and the category of coherent torsion
sheaves $\Coh^\infty(\EE)$.
\end{theorem}

This correspondence between torsion sheaves and semi-stable coherent sheaves can be
generalized to a relative setting of an elliptic fibration $\EE \lar \mathbb{T}$.
In  \cite{FMW} it was used
to construct vector bundles on
$\EE$ which are
semi-stable of degree zero on each fiber, see also \cite{BurbanKreussler2}.

As was shown in \cite{BurbanKreussler1},
the functor $\Phi$ is the trace of a certain  exact auto-equivalence
 of the derived category
$D^b(\Coh(\EE))$.
Using this equivalence of categories and a concrete description of
$\kk\llbracket x,y\rrbracket  /xy$--modules in terms of their projective resolutions,  
one can get a description
of semi-stable torsion
free sheaves of degree zero on a nodal Weierstra\ss{} curve in terms of \'etale coverings
\cite{FMmin, BurbanKreussler1}. In Section \ref{subsec:FMT}  we give a short
overview of some related
results without going into details.

\medskip
\medskip
\noindent
\textbf{Acknowledgement}.
Work on this article  has been supported by the DFG-Schwer\-punkt
"Globale Methoden in der komplexen Geometrie".
The second-named author would like to thank Bernd Kreu\ss{}ler for
numerous discussions of the results of this paper.

%\clearpage

\section{Vector bundles on smooth projective curves}\label{sec:smooth}

In this section we review some classical results about vector bundles on smooth
 curves. However, we provide non-classical proofs which, as we think, are simpler
and fit well in our approach to coherent sheaves over singular
curves.   The behavior  of the category of vector
bundles on a smooth projective curve $\XX$
is  controlled  by its genus  $g(\XX)$.

If $g(\XX)  = 0$ then $\XX$ is a projective line $\PP^1$ and any locally free sheaf
on it splits into
a direct  sum of line bundles $\kO_{\PP^1}(n)$, $n \in \mathbb{Z}$.
 This result, usually
attributed to Grothendieck \cite{Grothendieck}, was already known
in an equivalent  form
to
Birkhoff \cite{Birkhoff}. We found it quite instructive
to include   Birkhoff's algorithmic proof in our  survey.

A classification of vector bundles in the
case of smooth elliptic curves,  i.e. for 
$g(\XX) = 1$ was obtained  by Atiyah \cite{Atiyah}.
He has shown that the category of vector
bundles   on an elliptic curve $\XX$ is tame  and an indecomposable vector bundle
$\kE$ is
uniquely
determined by its rank $r$, its degree $d$ and a point of the curve $x \in \XX$. A modern,
and
%my change: from our point
in our opinion more conceptual way to prove  Atiah's result  uses the language of derived
categories and is due to Lenzing and Meltzer \cite{LenzingMeltzer}.
In the case of an  algebraically closed  field of
characteristics zero  an alternative description of indecomposable vector bundles via
\'etale coverings  was found by  Oda \cite{Oda}.
This classification was a cornerstone
in the proof of Polishchuk and Zaslow \cite{PolishchukZaslow}
of  Kontsevich's homological mirror symmetry conjecture for elliptic curves,
see also
\cite{Kreussler}.

The case of curves of genus bigger than one  is not considered in this survey. In this
situation
even the   category of semi-stable vector bundles of slope one is representation wild
and
the main attention is drawn to the study of various moduli problems
and properties of \emph{stable} vector bundles, see for example \cite{LePotier}.

\noindent
Throughout this section we do not require  any assumptions about
 the base field $\kk$.

\subsection{Vector bundles on the projective line}\label{subsec:projline}

We are going to  prove the following classical theorem.

\begin{theorem}[Birkhoff-Grothendieck]
Any vector bundle $\kE$ on the projective line $\PP^1$ splits into a direct sum of
line bundles:
$$
\kE \cong \bigoplus\limits_{n\in \mathbb{Z}} \kO_{\PP^1}(n)^{r_n}.
$$
\end{theorem}

A proof of this result  based on  Serre duality and vanishing theorems can be found in a book
of Le Potier \cite{LePotier}.
However, it is  quite interesting to give another, completely elementary
proof, based on a  lemma proven  by Birkhoff in 1913.

A projective line $\PP^1$ is a union of two affine lines $\mathbb{A}^1_i\ (i=0,1)$.  If
$(x_0:x_1)$ are homogeneous coordinates
 in $\PP^1$ then $\mathbb{A}^1_i=\{(x_0:x_1)| x_i\ne0\}$. The affine coordinate on
$\mathbb{A}^1_0$ is
$z =x_1/x_0$
 and on $\mathbb{A}^1_1$ it is  $z^{-1}=x_0/x_1$. Thus  we can identify
$\mathbb{A}^1_0$ with $\Spec(\kk[z])$ and
$\mathbb{A}^1_1$
 with $\Spec(\kk[z^{-1}])$,
 their intersection  is then $\Spec(\kk[z,z^{-1}])$. Certainly, any projective
module over
 $\kk[z]$ is free, i.e. all vector bundles over an affine line are trivial. Therefore to define a
vector bundle over $\PP^1$ one only has
 to prescribe its rank $r$ and a \emph{gluing matrix} $M\in GL(r, \kk[z, z^{-1}])$. Changing bases
in free modules over
 $\kk[z]$ and $\kk[z^{-1}]$ corresponds to the transformations $M\mapsto  T^{-1}MS$, where $S$ and $T$ are
invertible matrices
 of the same size,  over $\kk[z]$ and $\kk[z^{-1}]$ respectively.

\begin{proposition}[Birkhoff  \cite{Birkhoff}]
  For any matrix $M\in GL(r, \kk[z,z^{-1}])$ there are matrices
$S\in GL(r, \kk[z])$ and $T\in GL(r, \kk[z^{-1}])$ such
 that $T^{-1}MS$ is a diagonal matrix  $\mathsf{diag}(z^{d_1},\dots,z^{d_r})$.
 \end{proposition}

\noindent
\emph{Proof}. One can diagonalize the matrix $M$ in three steps.

\medskip
\noindent
Step 1.  Reduce the matrix $M = (a_{ij})$ to a  lower triangular form with diagonal entries
$a_{ii} = z^{m_i}$, where $m_i \in \mathbb{Z}$ and  $m_1 \le  m_2 \le \dots \le m_r$.
Indeed, since $\kk[z]$ is a discrete valuation ring, using invertible transformations of columns
 over  $\kk[z]$   we can reduce the first row $(a_{11}, a_{12}, \dots, a_{1r})$ of  $M$
to the form   $(a_1, 0, \dots, 0)$, where $a_1$ is the greatest common divisor of
$a_{11}, a_{12}, \dots, a_{1r}$.

Let $M'$ be the $(r-1) \times (r-1)$  matrix formed by the entries
$a_{ij}$, $i,j \ge 2$. Since  $\det(M) = a_1 \cdot \det(M')$ and $\det(M)$ is a unit in
$\kk[z, z^{-1}]$, it implies that  $a_1 = z^{m_1}$ and $\det(M')$ is a unit
in $\kk[z, z^{-1}]$ too.
Then we proceed with the matrix $M'$ inductively.
Note, that the diagonal entries  can be always reordered  to satisfy   $m_1 \le  m_2 \le \dots \le m_r$.

\medskip
\noindent
Step 2.  Consider the case of  a  lower-triangular $(2\times 2)$-matrix
$$M =
\left(
\begin{array}{cc}
z^{m} & 0 \\
p(z,z^{-1}) & z^{n}
\end{array}
\right)
$$
with  $m \le n$. We show by induction on the difference
$n - m$ that $M$ can be diagonalized
performing  invertible transformations of rows over $\kk[z^{-1}]$ and invertible
transformations of columns  over $\kk[z]$.

If $m = n$ then we can simply kill the entry $p = p(z, z^{-1})$. Assume now that
 $m < n$ and $p  \ne 0$.
Without loss of generality we may suppose  that $p \in \langle z^{m+1},\dots, z^{n-1}\rangle$.
Therefore there exist two mutually prime polynomials $a$ and $b$ in $\kk [z]  $ such that
$a p + b z^n = z^d$ and $m < d < n$. Then $\left(
\begin{array}{cc}
a & z^{n-d} \\
b & p/z^d
\end{array}
\right)$
belongs to $GL(2, \kk[z])$ and
$$
\left(
\begin{array}{cc}
z^{m} & 0 \\
p  & z^{n}
\end{array}
\right)
\left(
\begin{array}{cc}
a & z^{n-d} \\
b & p/z^d
\end{array}
\right) =
\left(
\begin{array}{cc}
z^{m} a & z^{n+m-d} \\
z^d  &  0
\end{array}
\right).
$$
In order
to conclude the induction step it remains to note that $|n + m -2d| < |n - m|$.
\medskip

\noindent
Step 3. Let $M$ be a lower-diagonal matrix with the diagonal elements\\
$z^{m_1},\dots, z^{m_r}$ with  $m_1 \le m_2 \le \dots \le m_r$.
We show by induction on
$\sum\limits_{i,j =1}^r  |m_i - m_j|$ that   $M$ can be diagonalized.
This  statement is obvious  for $\sum\limits_{i,j =1}^r|m_i - m_j| = 0.$
Assume that  $\sum\limits_{i,j =1}^r |m_i - m_j| = N  > 0$.
Introduce an ordering
on the set $\{(i,j)|1 \le j \le i \le r\}$:
 $$(2,1) < (2,3) < \dots < (r-1,r) < (3,1) < \dots < (r-2, r) < \dots < (r,1).$$
Let $(i_0,j_0)$ be the smallest pair such that $a_{i_0 j_0} \ne 0$. Then we can apply
the algorithm  from the Step 2 to the $(2 \times 2)$ matrix formed by the entries
$(j_0,j_0), (j_0,i_0),\\ (i_0,j_0)$ and $(i_0,i_0)$ to diminish the sum
$\sum\limits_{i,j =1}^r |m_i - m_j|$.
This  completes the proof of Birkhoff's lemma.

Now it remains to note that
$1\times 1$ matrix $(t^d)$ defines the line bundle $\kO_{\PP^1}(-d)$.
This  implies   the statement about   the splitting of a vector bundle on a projective line
into a direct sum of line bundles.

\subsection{Projective curves of arithmetic genus bigger then one
are vector bundle wild}\label{subsec:wild}

In this subsection we are going to prove the following

\begin{theorem}[see \cite{DrozdGreuel}]
Let $\XX$ be an irreducible projective curve of arithmetic genus $g(\XX) > 1$ over an
algebraically closed field $\kk$.
Then the
abelian category of semi-stable vector bundles of slope\footnote{The slope
of a coherent sheaf $\kF$ is $\mu(\kF) = \frac{deg(\kF)}{rk(\kF)}$.} one is representation wild.
\end{theorem}

In order to show the wildness of a category $\mathsf{A}$ one frequently uses the following
lemma:

\begin{lemma}
Let $\mathsf{A}$ be an abelian category, $M, N \in {\rm Ob}(\mathsf{A})$ with
$\Hom_{\mathsf{A}}(M,N)  = 0$ and
 $\xi$, $\xi' \in \Ext^{1}_{\mathsf{A}}(N,M) $  two extensions
$$
 \xi : 0 \lar M \stackrel{\alpha}\lar K \stackrel{\beta}\lar N \lar 0,
$$
$$
 \xi' : 0 \lar M \stackrel{\alpha'}\lar K' \stackrel{\beta'}\lar N \lar 0.
$$
Then $K \cong K'$ if and only if there exist  two isomorphisms
$ f : M\lar M$ and $g: N \lar N$ such that  $f\xi = \xi' g$.
\end{lemma}

\noindent
\emph{Proof}.
The statement is clear in one direction: if
$f\xi = \xi' g$ in $\Ext^{1}_{\mathsf{A}}(N,M)$, then   $K \cong K'$ by the $5$-Lemma.

Now suppose $K \cong K'$ and let  $h: K \lar K'$ be an isomorphism.
Then  $im(h\alpha)$ is a subobject of $im(\alpha')$. Indeed, otherwise the map
$\beta' h \alpha: M\lar N$ would be  non-zero, a contradiction.
Therefore we get the following commutative
diagram:

$$
%\begin{tabular}{p{2.1cm}c}
%&
\xymatrix
{
 0 \ar[r] & M\ar[r]^\alpha \ar[d]^f & K \ar[d]^{h} \ar[r]^\beta & N \ar[r] \ar[d]^{g} & 0 \\
0 \ar[r]& M \ar[r]^{\alpha'}  & K' \ar[r]^{\beta'} & N \ar[r]  & 0. \\
}
$$
%\end{tabular}

In the same way we proceed with  $h^{-1}$. Hence $f$ is an isomorphism, what proves the lemma.

\medskip

Let us come back to the proof of the theorem.
Suppose now that $\XX$ is an irreducible projective curve of arithmetic genus $g>1$, $
\kO:=\kO_{\XX}$.
Then for any  two points
 $x\ne y$ from $\XX$ we have $\Hom(\kO(x),\kO(y))\cong
 H^0(\XX,\kO(y-x))=0$ and the Riemann-Roch theorem  implies that
 $\Ext^1(\kO(x),\kO(y))\cong  H^1(\XX,\kO(y-x))\cong \kk^{g-1}$. Fix 5 different points $x_1,\dots,x_5$ of
the curve $\XX$, choose
 non-zero elements $\xi_{ij}\in \Ext^1(\kO(x_j),\kO(x_i))$ for $i\ne j$
and consider vector bundles $\kF(A,B)$,
 where $A,B\in Mat(n\times  n, \kk)$  and $\kF(A,B)$ is given as an extension
 $$
   0\lar  \underbrace{ (\kO(x_1)\oplus \kO(x_2))^n}_{\kB} \lar \kF(A,B)\lar
   \underbrace{(\kO(x_3)\oplus \kO(x_4)\oplus \kO(x_5))^n}_{\kA} \lar 0
 $$
 corresponding to the element $\xi(A,B)$ of $\Ext^1(\kA,\kB)$ presented by the matrix
 $$
\left(
   \begin{array}{ccc}
 \xi_{13}I &\xi_{14}I &\xi_{15}I\\
\xi_{23}I&\xi_{24}A&\xi_{25}B
\end{array}
\right),
 $$
where $I$ denotes the unit $n\times  n$ matrix. If $(A',B')$ is another pair of matrices,
and $\kF(A,B)\to\kF(A',B')$ any morphism, then the previous lemma implies
 that there are morphisms $\phi:\kA\to\kA'$ and
$\psi:\kB\to\kB'$ such
 that $\psi\xi(A,B)=\xi(A',B')\phi.$
Now one can easily deduce that $\Phi= \mathsf{diag}(S,S,S)$ and $\Psi=\mathsf{diag}(S,S)$
for some  matrix $S\in Mat(n \times  n, \kk)$ such that $SA=A'S$ and $SB=B'S$.

 If we consider a pair of matrices  $(A,B)$ as a representation of the free algebra
$\kk\langle x,y\rangle$ in 2
generators, the correspondence $(A,B)\mapsto
 \kF(A,B)$ becomes a \emph{full, faithful and  exact} functor
$\kk\langle x,y\rangle-\mod  \lar \mathsf{VB}(\XX)$. In particular, it maps non-isomorphic modules
 to non-isomorphic vector bundles and indecomposable modules to indecomposable vector bundles.
Using the terminology of
 representation theory of algebras, we say in this situation
that the curve $\XX$ is \emph{vector bundle wild}.
For a precise definition of wildness we refer to \cite{DrozdGreuel}.

 Recall that the algebra $\kk\langle x,y\rangle$ here can be replaced by
\emph{any} finitely generated algebra
$\Lambda = \kk\langle a_1,\dots,a_n\rangle$. Indeed,
 any $\Lambda$--module $M$ such that $dim_{\kk}(M)=m$ is given by a set of matrices
$A_1,\dots, A_n$  of size $m\times  m$.
 One gets a full, faithful and  exact functor $\Lambda-\bmod \lar \kk\langle x,y\rangle-\bmod$ mapping
the module $M$ to the
$\kk\langle x,y\rangle$--module of
 dimension $m\cdot n$ defined by the pair of matrices
 $$X=
\left(
\begin{array}{cccc}
   \lambda_1 I & 0 & \dots & 0\\
    0&\lambda_2 I&  \dots   &0      \\
    \vdots    &\vdots & \ddots & \vdots       \\
    0 &  0&  \dots & \lambda_n I
\end{array}
\right)
\quad
Y=
\left(
\begin{array}{cccc}
 A_1 &I     &\dots &0\\
 0   &A_2 & \dots&0 \\
\vdots    &\vdots  & \ddots &   \vdots  \\
0& 0&  \dots &  A_n
\end{array}
\right),
 $$
 where $\lambda_1, \dots, \lambda_n$ are different elements of
 the field $\kk$. Thus a classification
of vector bundles
over $\XX$ would imply a
 classifications of \emph{all} representations of \emph{all} finitely generated algebras,
a goal that perhaps nobody considers
 as achievable (whence the name ``wild'').

\subsection{Vector bundles on elliptic curves}\label{subsec:ellcurves}

In this subsection we shall discuss a classification
of indecomposable coherent sheaves  over smooth elliptic curves. Modulo some facts
about derived categories we give a self-contained proof of Atiyah's classification
of indecomposable vector bundles which is probably simpler than the original one.

\begin{definition}
An elliptic curve $\EE$ over a field $\kk$ is a smooth projective curve of genus one
having a $\kk$--rational point $p_0$.
\end{definition}

The category  $\Coh(\EE)$  of coherent sheaves on an elliptic curve $\EE$
has the following  properties, sometimes called
\emph{``the dimension one Calabi-Yau property''}:

\begin{itemize}
\item It is abelian, $\kk$--linear, $\Hom$-finite,
noetherian and of  global dimension one.
\item Serre Duality: for any two coherent sheaves $\kF$ and $\kG$ on $\EE$ there is an
isomorphism
$$
\Hom(\kF, \kG) \cong \Ext^1(\kG, \kF)^*,
$$
functorial in both arguments.
\end{itemize}

\noindent

\noindent
It is interesting to note that these properties almost  characterize the category of
coherent sheaves   on an elliptic curve:

\begin{theorem}[Reiten\,--\,van den Bergh \cite{ReitenvandenBergh}]
Let $\kk$ be an algebraically closed field and
$\sf{A}$  an indecomposable abelian Calabi-Yau category of dimension one.
Then $\sf{A}$ is equivalent either to the category of finite-dimensional
$\kk\llbracket t\rrbracket  $--modules or   to the category of coherent sheaves on  an
elliptic curve $\EE$.
\end{theorem}

This theorem characterizes Calabi-Yau abelian categories of global dimension one.
We shall need one more formula  to proceed with a  classification of
indecomposable coherent sheaves.

\begin{theorem}[Riemann--Roch formula]
\label{RR}
For any two coherent sheaves $\kF$ and $\kG$ on an elliptic curve
$\EE$ there is an integral bilinear Euler form
$$
\langle \kF, \kG\rangle: = dim_{\kk} \Hom(\kF, \kG) - dim_{\kk} \Ext^1(\kF, \kG) =
\left|
\begin{array}{cc}
deg(\kG) & deg(\kF) \\
rk(\kG)  & rk(\kF)
\end{array}
\right|.
$$
In particular, $\langle\,\, , \, \rangle$  is anti-symmetric:
 $\langle \kF, \kG\rangle = - \langle\kG, \kF\rangle$.
\end{theorem}

%\noindent
Now we  are ready to
 start with the  classification of indecomposable coherent sheaves.

\begin{theorem}[Atiyah]
Let $\EE$ be an elliptic curve over  a field $\kk$. Then
\begin{enumerate}
\item Any indecomposable coherent sheaf $\kF$ on $\EE$ is semi-stable.
\item If $\kF$ is semi-stable and indecomposable then all its Jordan-H\"older
factors are isomorphic.
\item A coherent sheaf $\kF$ is stable if and only if $\End(\kF) = \boldsymbol{K}$, where
$\kk \subset \boldsymbol{K}$ is some  finite field extension.
\end{enumerate}
\end{theorem}

%\noindent
\begin{proof} It is well-known that any coherent sheaf $\mathcal{F}\in\Coh(E)$
has a Harder-Narasim\-han filtration
$$0\subset \mathcal{F}_{n} \subset \ldots \subset \mathcal{F}_{1} \subset
\mathcal{F}_{0} = \mathcal{F}$$
whose factors $\mathcal{A}_{\nu} := \mathcal{F}_{\nu}/\mathcal{F}_{\nu+1}$ are
semi-stable with decreasing slopes $\mu(\mathcal{A}_{n})>
\mu(\mathcal{A}_{n-1}) > \ldots > \mu(\mathcal{A}_{0})$.
Using the definition of semi-stability, this implies\\
$\Hom(\mathcal{A}_{\nu+i}, \mathcal{A}_{\nu})
= 0$ for all $\nu\ge 0$ and $i>0$. Therefore,
$$\Ext^{1}(\mathcal{A}_{0},\mathcal{F}_{1}) \cong
\Hom(\mathcal{F}_{1}, \mathcal{A}_{0})^{\ast} =0,$$ and the exact sequence
$0\rightarrow \mathcal{F}_{1} \rightarrow \mathcal{F} \rightarrow
\mathcal{A}_{0} \rightarrow 0$ must split.
In particular, if $\mathcal{F}$ is indecomposable, we have $\mathcal{F}_{1}=0$
and $\mathcal{F}\cong \mathcal{A}_{0}$ and $\mathcal{F}$  is semi-stable.

The full sub-category of $\Coh(\EE)$ whose objects are the
semi-stable sheaves of a fixed slope is an abelian category in which any object
has a Jordan-H\"older filtration with stable factors.
If $\mathcal{F}$ and $\mathcal{G}$ are non-isomorphic stable sheaves which
have the same slope then $\Ext^1(\mathcal{F},\mathcal{G})= 0$. Based on this
fact  we  deduce that an indecomposable
semi-stable sheaf has all its Jordan-H\"older factors isomorphic to each
other.

It is well-known that any non-zero automorphism of
 a stable coherent sheaf $\mathcal{F}$ is invertible,
i.e. $\End(\mathcal{F})$ is a field $\boldsymbol{K}$.
Since $\EE$ is projective, the field extension
$\kk \subset \boldsymbol{K}$ is finite.
On a smooth
elliptic curve, the converse is true as well, which equips us with a useful
homological characterization of stability.

To see this,  suppose that
all endomorphism of $\mathcal{F}$  are invertible but $\kF$ is not stable.
This implies the existence of an epimorphism $\mathcal{F}\rightarrow \mathcal{G}$
with $\mathcal{G}$ stable and
$\mu(\mathcal{F})\ge \mu(\mathcal{G})$. Serre duality implies
$dim_{\kk} \Ext^{1}(\mathcal{G},\mathcal{F}) =
dim_{\kk} \Hom(\mathcal{F},\mathcal{G}) >
0$, hence, $\langle\mathcal{G},\mathcal{F}\rangle  =
dim_{\kk} \Hom(\mathcal{G},\mathcal{F}) -
dim_{\kk} \Ext^{1}(\mathcal{G},\mathcal{F}) < dim_{\kk}
\Hom(\mathcal{G},\mathcal{F})$. By Riemann-Roch formula
$\langle \mathcal{G},\mathcal{F}\rangle = \big(\mu(\mathcal{F}) -
\mu(\mathcal{G})\big)\rk(\mathcal{F})\rk(\mathcal{G}) > 0$, thus
$\Hom(\mathcal{G},\mathcal{F})\ne 0$. But this produces a
non-zero composition $\mathcal{F}\rightarrow \mathcal{G} \rightarrow
\mathcal{F}$ which is not an isomorphism, in contradiction to the assumption
that $\End(\kF)$ is a field.
\end{proof}

\begin{remark}
Usually one speaks about stability of vector bundles on projective varieties
in the case of an algebraically closed field of characteristics zero.
However, due to a result of Rudakov
\cite{Rudakov} one can introduce a stability notion for fairly general abelian
categories.
\end{remark}

The following classical fact was, probably first, proven  by Dold \cite{Dold}:

\begin{proposition}
Let $\mathsf{A}$ be an abelian category of global dimension one
 and $\kF$ an object of the derived category
$D^b(\mathsf{A})$. Then there is an isomorphism
$\kF \cong \bigoplus\limits_{i\in \mathbb{Z}} H^i(\kF)[-i]$,
i.e. any object of $D^b(\mathsf{A})$  splits into
a direct sum of its homologies.
\end{proposition}

This proposition in particular means that the derived category $D^b(\mathsf{A})$ of a
hereditary abelian  category $\mathsf{A}$ and the abelian category $\mathsf{A}$
itself have the same representation type.
However, it turns out that  the derived category has a richer structure and more symmetry
then  the corresponding abelian category.

First of all note that the group  $\Aut(D^b(\Coh(\EE)))$  acts on the $K$--group
 $K(\EE)$ of $\Coh(\EE)$ preserving  the Euler form $\langle\,\, ,\, \rangle$. Hence,  it leaves
invariant
the radical  of the Euler form
$\mathsf{rad}\langle\,\, , \,  \rangle = \{\kF \in K(\EE)| \langle\kF, \,\, \rangle = 0\}$
and induces an action on $K(\EE)/\mathsf{rad}\langle \,\,,\,\,\rangle$.
Since by Riemann-Roch theorem $Z: K(\EE)/\mathsf{rad}\langle\, ,\, \rangle \stackrel{Z}\lar
 \mathbb{Z}^2$ is an isomorphism, where $Z(\kF) := (\rk(\kF), \deg(\kF)) \in \mathbb{Z}^2$,
we get a group homomorphism $\Aut(D^b(\Coh(\EE))) \lar SL(2,\mathbb{Z})$. We call the
 pair $Z(\kF) \in \mathbb{Z}^2$  the \emph{charge} of $\kF$.

\begin{theorem}[Mukai, \cite{Mukai}]
Let $\EE$ be an elliptic curve. Then the  group homomorphism
$\Aut(D^b(\Coh(\EE))) \lar SL(2,\mathbb{Z})$ is surjective.
\end{theorem}

\noindent
\emph{Proof}. By the definition of an elliptic curve
there is a $\kk$--rational point $p_0$ on $\EE$ inducing
an exact equivalence $\kO(p_0)\otimes -$.
Let $\kP = \kO_{\EE\times \EE}\bigl(\Delta- (p_0\times\EE)-(\EE\times p_0)\bigr)$
then the Fourier-Mukai transform
$$\Phi_{\mathcal{P}}:D^{b}(\Coh(\EE)) \rightarrow D^{b}(\Coh(\EE)),
 \quad \Phi_{\mathcal{P}}(\kF) =
       \boldsymbol{R}\pi_{2\ast}(\mathcal{P}\otimes  \pi^{\ast}_{1}\kF)$$
is an exact auto-equivalence of the derived category, see \cite{Mukai}.
The actions of $\kO(p_0)\otimes -$ and $\Phi_{\mathcal{P}}$ on
$K(\EE)/\mathsf{rad}\langle\,\, ,\, \rangle$ in the basis $\{[\kO],[\kk(p_0)]\}$
are  given by the matrices
%my change: array
$
\left(
\begin{smallmatrix}
1 & 0 \\
1 & 1
\end{smallmatrix}
\right)
$
and $
\left(
\begin{smallmatrix}
0 & 1 \\
-1 & 0
\end{smallmatrix}
\right
)
$
which are known to generate $SL(2,\mathbb{Z})$. This shows
the claim.

\medskip
The technique of derived categories makes it easy to give a classification
of indecomposable coherent sheaves on an elliptic curve.

\begin{theorem}\label{thm:cohshcl}
Let $\kF$ be an indecomposable coherent sheaf on an elliptic curve $\EE$.
Then there exits a torsion sheaf $\kT$ and an exact auto-equivalence
$$\Phi \in \Aut(D^b(\Coh(\EE)))~~ such~ that~~\kF \cong \Phi(\kT).$$
\end{theorem}

\noindent
\emph{Proof}.
Let $\kF$ be an indecomposable coherent sheaf on $\EE$ with the charge  $Z(\kF) = (r, d)$,
$r > 0$.
Let $h = g.c.d.(r,d)$ be
%my change:
the greatest common divisor,
then there exists a matrix $F \in SL(2,\mathbb{Z})$ such that
$F \left(\begin{smallmatrix} r \\ d\end{smallmatrix}\right)  = \left(\begin{smallmatrix} 0 \\
h\end{smallmatrix}\right)$.
We can lift the matrix $F$ to an auto-equivalence $\Phi \in \Aut(D^b(\Coh(\EE))$, then
$Z(\Phi(\kF)) = \left(\begin{smallmatrix} 0 \\
h\end{smallmatrix}\right)$.
Since
$\Aut(D^b(\Coh(\EE)))$ maps indecomposable objects of the derived category to indecomposable
ones  and  since the only indecomposable objects in the derived category
 are shifts of indecomposable coherent sheaves,
we can conclude  that $\Phi(\kF)$ is isomorphic to
a shift of some indecomposable sheaf of rank zero,  what proves the theorem.

\medskip
Let $M_{\EE}(r,d)$ denote  the set of indecomposable vector bundles on $\EE$ of rank
$r$ and degree $d$.

\begin{theorem}[Atiyah]
Let $\EE$ be an elliptic curve. Then for any integer $h > 0$ there exists
a unique indecomposable
vector bundle $\kF_h \in M_{\EE}(h,0)$ such that $H^0(\kF_h) \ne 0$. The vector
bundles $\kF_h$
are  called \emph{unipotent}. Moreover, the following properties hold:
\begin{enumerate}
\item $H^0(\kF_h) = H^1(\kF_h) = \kk$ for all $h \ge 1$.
\item If $char(\kk) = 0$ then $\kF_h \cong Sym^{h-1}(\kF_2)$.  Moreover
$$
\kF_e \otimes \kF_f \cong \bigoplus\limits_{i= 0}^{f-1} \kF_{e+ f- 2i -1}
$$
\end{enumerate}
\end{theorem}

\noindent
\emph{Sketch of the proof}.
Since  $\kF_h$ is  indecomposable  of degree zero, it
has a unique
Jordan-H\"older factor $\kL \in \Pic^0(\EE)$. From the assumption $\Hom(\kO, \kF_h) \ne 0$
we conclude that
$\kL \cong \kO$, so  each
bundle  $\kF_h$ can be obtained by recursive  self-extensions
of the structure sheaf.
Since by Theorem \ref{thm:cohshcl} the category of semi-stable vector bundles of degree zero
is equivalent to the category of
torsion sheaves, we conclude that the category of semi-stable sheaves with the
Jordan-H\"older factor $\kO$ is equivalent to the category of finite-dimensional
$\kk\llbracket t\rrbracket  $-modules. The
exact sequence
$$
0 \lar \kF_{h-1} \lar \kF_h \lar \kO \lar 0.
$$
corresponds via   Fourier-Mukai transform $\Phi_\kP$ to
$$
0 \lar \kk\llbracket t\rrbracket  /t^{h-1} \lar \kk\llbracket t\rrbracket  /t^{h} \lar \kk \lar 0.
$$
In the same way we  conclude that $\Hom(\kO, \kF_h) =
\Hom_{\kk\llbracket t\rrbracket  }(\kk,  \kk\llbracket t\rrbracket  /t^{h}) = \kk$.
Moreover, one can show that
$\Phi_\kP(\kk\llbracket t\rrbracket  /t^{e} \otimes_{\kk} \kk\llbracket t\rrbracket  /t^{f}) \cong \kF_e \otimes \kF_f$,
hence
we have the same rules for the decomposition of the tensor product of  unipotent vector
bundles and of  nilpotent Jordan cells, see \cite{Atiyah, Oda, Mukai, PloogHein}.

\medskip
\begin{remark}
Atiyah's original proof from 1957 was written at the time when the formalism of derived
and triangulated categories was not developed  yet. However, his construction of a bijection
between
$M_{\EE}(r,d)$ and $M_{\EE}(h,0)$ corresponds exactly to the action of the group of
exact auto-equivalences
of the derived category of coherent sheaves
 on the set of indecomposable objects. This was probably for the first time
observed   by Lenzing and Meltzer in   \cite{LenzingMeltzer}. For  further
elaborations,  see  \cite{PolishchukBook, PloogHein, BurbanKreussler3}.
\end{remark}

Actually, Atiyah's description of indecomposable vector bundles on an elliptic curve
$\EE$ is  more precise.

\begin{theorem}[Atiyah]\label{thm:Atiyah}
Let $\EE$ be an elliptic curve over an algebraically closed field $\kk$. For
any pair of coprime integers  $(r,d)$ with $r > 0$ pick up  some
$\kE(r,d) \in M_{\EE}(r,d)$. Then
\begin{enumerate}
\item $M_{\EE}(r,d) = \{\kE(r,d) \otimes \kL| \, \kL \in \Pic^0(\EE)\}$.
\item $\kE(r,d) \otimes \kL \cong \kE(r,d)$ if and only if $\kL^r \cong \kO$.
\item The map $det: M_{\EE}(r,d) \lar M_{\EE}(1,d)$ is a bijection.
\item If  $char(\kk) = 0$,  then $\kF_h\otimes - : M_{\EE}(r,d) \lar M_{\EE}(rh,dh)$
      is a bijection.
\end{enumerate}
\end{theorem}

\begin{remark}
If $\kk$ is not algebraically closed and $\EE$
is a smooth projective curve of genus one over $\kk,$  without $\kk$--rational points, then
 we miss the generator $\kO(p_0)$ in the
group of exact auto-equivalences $D^b(\Coh(X))$ and  the method used for elliptic curves
can not be immediately applied.
This problem was solved in a paper of Pumpl\"un \cite{Pumpluen}.
\end{remark}

We may sum up  the discussed properties  of indecomposable coherent sheaves on
elliptic curves:

\begin{proposition}
\label{propCoh}
Let $\EE$ be an elliptic curve over a field $\kk$.  Then
\begin{enumerate}
\item Any indecomposable coherent sheaf $\kF$ on $\EE$ is semi-stable with a unique stable
Jordan-H\"older factor.

\item An indecomposable vector bundle
is determined by its charge $(r, d)\in \mathbb{Z}^2$ and a closed point
$x$ of the curve $\EE$.

\item Let $\Coh^\nu(\EE)$ be the category of
semi-stable sheaves of slope $\nu$. Then
 for any $\mu, \nu \in \mathbb{Q} \cup\{\infty\}$ the
abelian categories  $\Coh^\nu(\EE)$ and $\Coh^\mu(\EE)$ are equivalent and this equivalence is
induced by an auto-equivalence of $D^b(\Coh(\EE))$.
\item
In particular, each category
$\Coh^\mu(\EE)$ is equivalent to the category of coherent torsion sheaves.

\item If $\kF \in \Coh^\nu(\EE), \kG \in \Coh^\mu(\EE)$ and   $\nu < \mu$ then
$\Ext^1(\kF, \kG) = 0$ and $$dim_{\kk} \Hom(\kF, \kG) =
deg(\kG)rk(\kF) - deg(\kF)rk(\kG).$$
The case $\nu >\mu$ is dual by  Serre duality.
\end{enumerate}
\end{proposition}

This gives a pretty complete description of the category of coherent sheaves on an elliptic
curve and of its derived category. However, in applications one needs another
description of indecomposable vector bundles, see \cite{Polishchuk, PolishchukZaslow}.

The following form of Atiyah's
classification is due to Oda \cite{Oda}.
It  was used by Polishchuk and Zaslow in their proof of the
homological mirror symmetry  conjecture for elliptic curves, see
\cite{PolishchukZaslow, Kreussler}.

\begin{theorem}[Oda]
Let $\kk = \mathbb{C}$ and $\EE = \EE_\tau = \mathbb{C}/\langle1,\tau\rangle$
be an elliptic curve, $\kE \in M_ {\EE}(rh,dh)$ an indecomposable vector bundle, where
$g.c.d.(r,d) = 1$. Then there exists a unique
 line bundle $\kL$ of degree $d$ on $\EE_{r\tau}$ such that
$$\kE \cong \pi_{*}(\kL)\otimes \kF_h \cong \pi_{*}(\kL \otimes \kF_h),
$$
where $\pi: \EE_{r\tau} \lar \EE_\tau$ is an \'etale covering of degree $r$.
\end{theorem}

\noindent
\emph{Proof}. Let $\kL$ be a line bundle on $\EE_{r\tau}$ of degree $d$.
Since the morphism $\pi$ is \'etale, $\pi_{*}(\kL)$ is a vector bundle on $\EE_\tau$ of
rank $r$. The Todd class of an elliptic curve is trivial, hence
by Grothendieck-Riemann-Roch theorem we obtain
$deg(\pi_{*}(\kL)) = deg(\kL) = d$.

Now let us show  that $\pi_{*}(\kL)$ is indecomposable. To do this it suffices
to prove  that $\End(\pi_{*}\kL) = \mathbb{C}$.
Consider the fiber product diagram
$$
\begin{CD} \widetilde{\EE} @>p_{1}>> \EE_{r\tau} \\ @Vp_{2}VV
          @VV\pi_{1}V \\ ~\EE_{r\tau} @>\pi_{2}>> ~\EE_\tau
\end{CD}
$$
One can easily check that $\widetilde{\EE}$ is a union  of $r$ copies of the elliptic curve
$\EE_{r\tau}$: $\widetilde{\EE} = \coprod\limits_{i= 1}^r \EE_{r\tau}^i$, where
each $p_1^{(i)}:  \EE_{r\tau}^i  \lar
\EE_{r\tau}$, $i = 1,\dots,r$ can be chosen to be the identity
map and $p_2^{(i)}(z) = z + \frac{i}{r}\tau$.

Since all morphisms $\pi_i, p_i, i = 1,2$ are affine and flat,
the  functors $\pi_{i*}, \pi_{i}^*, p_{i*},
p_{i}^* $ are exact. Moreover,  $p_i^!  = p_i^*,$
since the canonical sheaf of an elliptic
curve is trivial.\footnote{Recall that
if $f:\XX \to  \YY$ is a finite morphism of Gorenstein projective schemes then
$\Hom_{\XX}(\kF,f^{!}\kG) \cong
            \Hom_{\EE}(f_{\ast}\kF,\kG)$
for any coherent sheaf $\kF$ on $\XX$ and a coherent sheaf  $\kG$ on $\YY.$
}
Using the base change isomorphism and Grothendieck duality we have
\begin{align*}
 \Hom_{\EE\tau}\bigl(\pi_{1*}\kL ,
\pi_{2*}\kL\bigr) \cong \Hom_{\EE_{r\tau}} \bigl(\pi_{2}^{*} \pi_{1*}\kL, \kL\bigr)\\
\cong  \Hom_{\EE_{r\tau}} \bigl(p_{2*} p_{1}^{*}
\kL, \kL\bigr) \cong \Hom_{\widetilde{\EE}}\bigl(p_{1}^{*} \kL,
p_{2}^{*}\kL\bigr).
\end{align*}
It remains to note that $\Hom_{\widetilde{\EE}}\bigl(p_{1}^{(i)*} \kL,
p_{2}^{(i)*}\kL\bigr) = 0$ for $i \ne 0$.

If $\kE$ is an  indecomposable vector bundle on $\EE_\tau$ of rank $rh$ and degree $dh$, then
by Theorem \ref{thm:Atiyah} there exists $\kM \in \Pic^0(\EE_\tau)$ such that
$\kE \cong \pi_{*}(\kL) \otimes \kM \otimes \kF_h$.
By the projection formula   $\kE \cong \pi_{*}(\kL  \otimes \pi^*\kM)\otimes \kF_h$.
Moreover,  passing to an \'etale covering kills the ambiguity in the choice of $\kM$.

It remains to show that $\pi^*(\kF_h) \cong \kF_h.$ To do this  it suffices
see that $\pi^*(\kF_2) \cong \kF_2,$ since $\kF_h \cong Sym^{h-1}(\kF_2)$ and
the inverse image commutes with all tensor operations.
The only property we have to check is that $\pi^*(\kF_2)$ does not split.
It is equivalent to say that the map $\pi^*: H^1(\kO_{\EE_\tau}) \lar
H^1(\kO_{\EE_{r\tau}})$ is non-zero.

%tut bude diagramma
%\marginpar{here will be a diagram }

%my change:
This follows from the  commutativity of the diagram:
$$
\newcommand{\A}
{\left(
\begin{smallmatrix}
1 & 0 \\
0 & r
\end{smallmatrix}
\right)
}
%{\left(\begin{smallmatrix }1&0\\ 0&r\end{smallmatrix}\right)}
\begin{xymatrix}
{
\mZ^2\ar@{=}[r]
 & H_1(\EE_\tau, \mZ)\ar[r]^{\cong} & H^1(\EE_\tau, \mZ) \ar@^{(->}[r]\ar[d]^{\pi^*} &
H^1(\EE_\tau, 0)\ar[d]^{\pi^*}\\
\mZ^2\ar@{=}[r]\ar[u]^{\A}
 & H_1(\EE_{r\tau}, \mZ)\ar[u]^{\pi_*}\ar[r]^{\cong}  & H^1(\EE_{r\tau}, \mZ)\ar@^{(->}[r]  &
H^1(\EE_{r\tau}, 0).\\
}
\end{xymatrix}
$$

%\clearpage
\section{Vector bundles and torsion free  sheaves on
singular curves of arithmetic genus one}\label{sec:singular}

In this paper we discuss
two approaches for  the study of the category of coherent
sheaves on a singular projective  curve of arithmetic genus one.
The first   uses the technique
of derived categories and
Fourier-Mukai transforms.  Its key point is that
any semi-stable torsion free sheaf on an irreducible Weierstra\ss{} curve can be
obtained from a torsion sheaf by applying an auto-equivalence of the derived category.
This technique can be generalized to the case of elliptic fibrations:
we can transform a family of torsion sheaves to a family of sheaves, which are
 semi-stable  on each fiber.

However, the approach via Fourier-Mukai transforms allows to describe only
semi-stable   sheaves.
In order to get a description of all  indecomposable  torsion free sheaves,  another
technique turns out to be useful. Namely, we relate vector bundles on a singular
rational curve $\XX$ and on its normalization $\widetilde{\XX} \stackrel{p}\lar \XX$.
The inverse image functor
%my change: VB
$p^*: \VB(\XX)\lar \VB(\widetilde{\XX})$ can map non-isomorphic bundles into
isomorphic ones.  The full information about the fibers of this map
is encoded in a certain \emph{matrix problem}.  In the case of an algebraically closed field
this approach leads to a very
concrete description of indecomposable vector bundles on cycles of projective lines
via \'etale coverings (no assumption on $char(\kk)$ is   needed).

Combining both  methods,   we get a quite complete description of the category
of torsion free sheaves  on  a nodal   Weierstra\ss{} curve.

\subsection{Vector bundles on singular curves via matrix problems}\label{subsec:mp}
\label{subsecTr}

Let  $\XX$ be a reduced projective curve over a field $\kk$. Introduce
the following notation:

\begin{itemize}
\item  $p: \widetilde{\XX} \lar \XX$  the normalization of $\XX$;
\item   $\kO=\kO_{\XX}$  and $\widetilde{\kO} = p_{*}\kO_{\widetilde{\XX}}$;
\item   $\,\kJ= Ann_{\kO}(\widetilde{\kO}/\kO)$ the conductor of $\kO$ in  $\widetilde{\mathcal O}$;
\item   $\,\kA=\kO/\kJ\,$  and $\,\widetilde{\kA}=\widetilde{\kO}/\kJ$.
\end{itemize}

\noindent
 Note that $\,\kA\,$ and $\,\widetilde{\kA}\,$ are skyscraper sheaves
supported at the singular locus of $\XX$.
Since the morphism $p$ is affine, $p_*$  identifies  the category of coherent sheaves
$\Coh(\widetilde{\XX})$ and the category $\Coh_{\widetilde{\kO}}$
of coherent modules on the ringed space $(\XX, \widetilde{\kO})$. Let $\ZZ$ be the subscheme of
$\XX$ defined by the conductor $\kJ$, $\widetilde\ZZ$ its scheme-theoretic pull-back on
$\widetilde\XX$ and $\kI = \kI_{\widetilde\ZZ}$ its ideal sheaf on $\widetilde\XX$.
Then $p_*$ also induces an equivalence
between  the category of $\kO_{\widetilde\XX}/\kI$--modules and the category of
$\widetilde\kA$--modules.

For a sheaf of algebras $\Lambda \in \{\kO, \widetilde\kO, \kA, \widetilde\kA \}$
on the topological space  $\XX$,
denote by $\TF_\Lambda$ the category of torsion free coherent $\Lambda$--modules and
by $\VB_\Lambda$ its full subcategory of locally free sheaves.
The usual way to deal with  vector bundles on a singular curve is to lift them
 to the normalization, and then work on a
smooth curve, see for example \cite{Seshadri}.
Passing to the normalization  we loose information about the isomorphism classes of objects
of $\VB_\kO$
 since non-isomorphic vector bundles
can have isomorphic inverse images.  In order to describe the fibers of the map
$\VB_{\kO} \lar \VB_{\widetilde\kO}$   and to be able to deal with  arbitrary torsion free sheaves
we introduce the  following definition:

\begin{definition}
\emph{The category of triples}   ${\mathsf T}_{\XX}$  is defined as follows:
\begin{enumerate}
\item Its objects are  triples
$(\widetilde{\kF}, {\kM}, \tilde{i})$, where $\widetilde{\mathcal F}$ is a
locally free  $\widetilde\kO$--module,
${\mathcal M}$  is a coherent   $\kA$--module and
$\tilde{i} : {\mathcal M}\otimes_{\kO} \widetilde{\kA} \lar
 \widetilde{\mathcal F} \otimes_{\widetilde\kO} \widetilde\kA $
is an epimorphism of $\widetilde\kA$--modules, which induces a monomorphism of $\kA$--modules
$
i: \kM \lar {\kM} \otimes_{\kO} \widetilde{\kA} \stackrel{\tilde{i}}\lar
\widetilde{\kF}\otimes_{\widetilde\kO} \widetilde\kA
$.

\item  A morphism
 $
(\widetilde{\kF}_{1}, {\kM_{1}}, \tilde{i}_{1})
\stackrel{(F, f)}\lar
(\widetilde{\kF}_{2}, {\kM_{1}}, \tilde{i}_{2})
 $
is given by a pair  $(F, f)$, where
$\widetilde{\kF}_{1} \stackrel{F}\lar \widetilde{\kF}_{2}$
is a morphism of $\widetilde{\kO}$-modules and
$
 {\kM_{1}} \stackrel{f}\lar
 {\kM_{2}}
$ is a morphism of  $\kA$-modules, such that the following diagram

\begin{tabular}{p{3cm}c}
 &
\xymatrix
{
{\mathcal M_{1}} \otimes_{\kO} \widetilde{\kA}  \ar[rr]^{\tilde{i}_1} \ar[d]_{\bar{f}} & &
\widetilde{\mathcal F}_{1} \otimes_{\widetilde \kO} \widetilde{\kA}
\ar[d]^{\bar{F}}  \\
{\mathcal M_{2}} \otimes_{\kO} \widetilde{\kA}  \ar[rr]^{\tilde{i}_2}  & &
\widetilde{\mathcal F}_{2} \otimes_{\widetilde \kO} \widetilde{\kA} \\
}
\end{tabular}

\noindent
is commutative in $\Coh_{\widetilde\kA}$, where $\bar{F} = F\otimes id$ and $\bar{f} =
\varphi\otimes id$.
\end{enumerate}
\label{dfTriples }
\end{definition}

\noindent
The main reason  to introduce the formalism of triples is the following theorem:

\begin{theorem}\label{thm:triples}
The functor
$
\TF_{\kO} \stackrel{\Psi} \lar {\mathsf T}_{\XX}
$
mapping a torsion free sheaf $\kF$ to the triple
$(\widetilde\kF, \kM, \tilde{i})$, where
$\widetilde{\kF} = \kF\otimes_\kO \widetilde{\kO}/\mathsf{tor}(\kF\otimes_\kO \widetilde{\kO})$,
$\kM = \widetilde{\kF}\otimes_{\kO} \kA$ and $
\tilde{i}: \kF \otimes_\kO \widetilde\kA  \lar \widetilde{\kF}\otimes_{\widetilde{\kO}} \widetilde\kA,
$
is an equivalence of categories. Moreover, the category of vector bundles
 $\VB_\kO$ is equivalent to
the full subcategory of ${\mathsf T}_{\XX}$ consisting
of those triples $(\widetilde\kF, \kM, \tilde{i})$, for which $\kM$ is a free
$\kA$--module and
$\tilde{i}$
is an isomorphism.
\end{theorem}

\noindent
\emph{Sketch of the proof}.
We construct the quasi-inverse functor
$
 {\mathsf  T}_{\XX}  \stackrel{\Psi'} \lar \TF_{\XX}
$
as follows.
Let
$(\widetilde{\mathcal F}, {\mathcal M}, \tilde{i})$ be some triple.
Consider  the pull-back diagram

\begin{tabular}{p{2.5cm}c}
&
\xymatrix
{0 \ar[r] & \kJ \widetilde{\kF} \ar[r] \ar[d]^{id} & \kF \ar[r] \ar[d] &
 \kM \ar[r] \ar[d]^{i} & 0 \\
0 \ar[r] & \kJ \widetilde{\kF} \ar[r]  & \widetilde{\kF} \ar[r]  &
 \widetilde{\kF}\otimes_{\widetilde{\kO}} \widetilde{\kA} \ar[r] & 0
}
\end{tabular}

\noindent
in the category of
 $\kO$--modules. Since the  pull-back is functorial, we get  a functor
${\mathsf T}_{\XX}  \stackrel{\Psi'} \lar \Coh(\XX)$. Since the map $i$ is injective,
$\kF \lar  \widetilde{\kF}$ is injective as well, so $\kF$ is torsion free.
It remains  to show that the functors  $\Psi$ and
$\Psi'$ are quasi-inverse to each other. We refer to \cite{DrozdGreuel} for the
details of the proof.

\begin{remark}\label{rem:triples}
There is a geometric way to interpret the above construction of the category of triples.
Let $\XX$ be a singular curve, $\widetilde\XX \stackrel{p}\lar \XX$ its normalization,
 $s: \mathbb{S} \lar \XX$  the inclusion of
the closed subscheme defined by
the conductor ideal and $\tilde{s}: \widetilde{\mathbb S} \lar \widetilde\XX$ its pull-back on the normalization. Consider the Cartesian  diagram
$$
\xymatrix{
\widetilde{\mathbb{S}} \ar[r]^{\tilde{s}} \ar[d]^{\tilde{p}}  & \widetilde\XX \ar[d]^{p}\\
\mathbb{S} \ar[r]^{s}  &  \XX.
}
$$
Theorem \ref{thm:triples}  says  that a torsion free sheaf $\kF$ on a singular  curve $\XX$
can be reconstructed from its ``normalization'' $p^*(\kF)/\mathsf{tor}(p^*\kF)$,
its pull-back $s^*\kF$ on $\mathbb{S}$  and the ``gluing map''
$\tilde{p}^* s^*\kF \lar \tilde{s}^* p^*\kF \lar
\tilde{s}^*\big(p^*\kF/\mathsf{tor}(p^*\kF)\big)$.
\end{remark}

\medskip
\medskip

Now let us see how this construction can be used to classify torsion free sheaves
on degenerations of elliptic curves.
Let $char(\kk) \ne 2$ and
$\EE$ be a nodal Weierstra\ss{} curve, given by the equation $zy^2 - x^3 -zx^2 = 0$,
$s = (0: 0: 1)$ its singular point, ${\PP}^1 = \widetilde{\EE} \stackrel{p}\lar \EE$ the
normalization map.
Choose coordinates on ${\PP}^1$ in such a way  that the preimages of $s$
are $0 = (0 : 1)$ and $\infty = (1 : 0)$.

The previous theorem says that a  torsion free sheaf
$\kF$ on the curve $\EE$ is uniquely determined by the corresponding triple
$\Psi(\kF) = (\widetilde{\kF}, \kM, \tilde{i})$.
Here $\widetilde{\kF}$ is a locally free
$\widetilde{\kO}$--module,
or as we have  seen, a  locally free $\kO_{\PP^{1}}$-module.
Using the notation $\widetilde{\kO}(n) = p_*(\kO_{\PP^1}(n))$,
due to  the  theorem of  Birkhoff-Grothendieck,
$\widetilde{\kF} \cong \bigoplus\limits_{n\in \mZ} \widetilde{\kO}(n)^{r_{n}}$.

 Since $\kA = \kO/\kJ = \kk(s)$ and $\widetilde{\kA} = \widetilde{\kO}/\kJ = (\kk\times\kk)(s)$,
the sheaf $\kM$   can be identified with its stalk at $s$ and the map
$\tilde{i} :  \kM\otimes_\kA  \widetilde\kA \lar    \widetilde\kF\otimes_{\widetilde\kO} \widetilde\kA$
can be viewed as a pair $(i(0), i(\infty))$ of
linear maps of $\kk$--vector spaces.
In order to write
$\tilde{i}$ in terms of matrices we  identify
$\widetilde\kO(n)\otimes_{\widetilde{\kO}}\widetilde\kA$ with
$p_*(\kO_{\PP^1}(n)\otimes_{\kO_{\PP^1}} \kO_{\PP^1}/\kI)$.  The choice of coordinates
on $\PP^1$ fixes two canonical sections $z_0$ and $z_1$ of $H^0(\kO_{\PP^1}(1))$
and we  use the trivializations $$\kO_{\PP^1}(n)\otimes \kI \lar \kk(0) \times  \kk(\infty)$$
 given
by
$\zeta\otimes 1 \mapsto (\zeta/z_1^n(0), \zeta/z_0^n(\infty))$.
Note, that this isomorphism  only  depends on the choice of coordinates of $\PP^1$.
In such a way we supply the $\widetilde\kA$--module
$\widetilde\kF\otimes \widetilde\kA = \widetilde\kF(0) \oplus \widetilde\kF(\infty)$ with a basis and get
isomorphisms
$\widetilde\kF(0) \cong \bigoplus\limits_{n\in \mZ} \kk(0)^{r_{n}}$
and
$\widetilde\kF(\infty) \cong \bigoplus\limits_{n\in \mZ} \kk(\infty)^{r_{n}}$.
With respect to  all  choices the morphism $\tilde{i}$
is given by two  matrices $i(0)$ and $i(\infty)$,
divided into horizontal blocks:

%my picture:
\begin{figure}[ht]
\hspace{2cm}
\begin{minipage}[b]{4cm}
\begin{align*}
\setlength\arrayrulewidth{0.3mm}
\begin{array}{ @{} c@{}     c  @{}c@{}  r }
\cline{1-2}
\fvr & %\hspace{3cm} d%\phantom{abrakadabradabra}
{~\mk{}\mk{} {\vdots}  \mk{}\mk{}~}
& \fvr & \\
%\cline{1-2}
%\fvr &  {\vdots} & \fvr\\
\cline{1-2}
\fvr &{\scriptstyle n-1} & \fvr & \\
\cline{1-2}
\fvr & {\scriptstyle n} & \fvr  &  \}{\scriptstyle r_n} \\
\cline{1-2}
\fvr &_{\scriptstyle n+1} & \fvr &\\
\cline{1-2}
\fvr & ^{\vdots} & \fvr &\\
\cline{1-2}
%\fvr & & \fvr\\
%\cline{1-2}
\end{array}
\end{align*}
\center {$i(0)$}
%\caption*{c)}
\end{minipage}
\begin{minipage}[b]{4cm}
\begin{align*}
\setlength\arrayrulewidth{0.3mm}
\begin{array}{ @{} c@{} c      @{}c@{} l }
\cline{1-2}
\fvr & %\hspace{3cm} d%\phantom{abrakadabradabra}
{~\mk{}\mk{}  {\vdots}  \mk{}\mk{}~}
& \fvr\\
%\cline{1-2}
%\fvr &  {\vdots} & \fvr\\
\cline{1-2}
\fvr &{\scriptstyle n-1} & \fvr\\
\cline{1-2}
\fvr & {\scriptstyle n} & \fvr & \} {\scriptstyle r_n} \\
\cline{1-2}
\fvr &_{\scriptstyle n+1} & \fvr\\
\cline{1-2}
\fvr & ^{\vdots} & \fvr\\
\cline{1-2}
%\fvr & & \fvr\\
%\cline{1-2}
\end{array}
\end{align*}
\center {$i(\infty)\phantom{BB}$}
%\caption*{c)}
\end{minipage}
\end{figure}

From the definition of the category of triples it follows that the matrices
$i(0)$ and $i(\infty)$ have to be of   full row rank and the transposed
matrix $(i(0)|i(\infty))^t$ has to be
monomorphic. Vector bundles on $\EE$ correspond to  invertible
square matrices $i(0)$ and $i(\infty)$.

Of course, for a fixed $\widetilde\kF = \bigoplus\limits_{n\in \mZ} \widetilde{\kO}(n)^{r_{n}}$ and
$\kM = \kk^N(s)$,  two different pairs of matrices $(i(0)|i(\infty))$ and $(i'(0)|i'(\infty))$
can define isomorphic torsion free sheaves on $\EE$. However, since the functor
$\Psi: \TF_\kO \lar {\mathsf T}_\XX$
preserves isomorphism classes of indecomposable objects,
two triples
$(\widetilde{\kF}, \kM, \tilde{i})$ and $(\widetilde{\kF}, \kM, \widetilde{i'})$
define isomorphic torsion free sheaves  if  and only if there are   automorphisms
$F : \widetilde{\kF} \lar \widetilde{\kF}$ and $f: \kM \lar \kM$
such that $\bar{F} \tilde{i} = \widetilde{i'}  \bar{f}$.

An endomorphism $F$ of
$\widetilde{\kF} = \bigoplus\limits_{n \in \mathbb{Z}} \widetilde{\kO}(n)^{r_n}$ can be written in a
matrix form: $F = (F_{kl})$, where $F_{kl}$ is
a $r_{l} \times r_{k}$--matrix with  coefficients in the vector space
$\Hom(\widetilde{\kO}(k), \widetilde{\kO}(l)) \cong
\kk[z_{0}, z_{1}]_{l-k}$. In particular, the matrix $F$
is lower triangular and the diagonal $r_{n}\times r_{n}$ blocks
$F_{nn}$ are just  matrices over $\kk$. The morphism
$F$ is an isomorphism if and only if
all $F_{nn}$ are invertible.
 Let $r = rank(\widetilde{\kF})$. With respect to the chosen trivialization
of $\kO_{\PP^1}(n)$ at $0$ and $\infty$ the map
$\bar{F} : \kk^{r}(0) \oplus \kk^{r}(\infty) \lar \kk^{r}(0) \oplus \kk^{r}(\infty)$
is given by the pair of matrices  $(F(0), F(\infty)) $
and  we have the following transformation rules
for the pair $(i(0), i(\infty))$:
$$
\bigl(i(0), i(\infty)\bigr) \mapsto
\bigl(F(0)^{-1} i(0) S, F(\infty)^{-1} i(\infty)S\bigr),
$$
where $F$ is an automorphism of
$\bigoplus\limits_{n \in \mathbb{Z}} \widetilde{\kO}(n)^{r_n}$ and
$S$ an automorphism of $\kk^N$.
Note, that the matrices
$F_{kl}(0)$ and $F_{kl}(\infty)$, $k,l \in \mZ, k > l$ can be \emph{arbitrary}
and
$F_{nn}(0) = F_{nn}(\infty)$ can be \emph{arbitrary invertible} for $n \in \mathbb{Z}$.
As a result we get   the following \emph{ matrix problem}.

\medskip
\medskip
\noindent
\textbf{Matrix problem for a nodal Weierstra\ss{} curve}.
We have two matrices $i(0)$ and $i(\infty)$ of the same size and both of full row rank.
Each of them is divided  into horizontal blocks labeled by  integers  (they
are called sometimes weights). Blocks of $i(0)$ and $i(\infty)$, labeled by the same integer,
have  the same size. We are allowed to   perform only
the following transformations:

\begin{enumerate}
\item We can simultaneously do any elementary transformations of columns of
$i(0)$
and $i(\infty)$.
\item We can simultaneously do any invertible elementary transformations of rows
inside of any two  conjugated horizontal blocks.
\item We can in each of the matrices $i(0)$ and
$i(\infty)$  independently add
a scalar multiple of any row with lower weight to any row with higher weight.
\end{enumerate}

The main idea is that we can transform the matrix $i$ into a
canonical form which is quite analogous to the Jordan normal form.

These types of matrix problems are well-known in  representation theory.
First they appeared in a work of Nazarova and Roiter \cite{NazRoi}
about the
classification of $\kk\llbracket x,y\rrbracket  /(xy)$--modules. They are called,
sometimes, ``Gelfand problems'' or ``representations of bunches of chains''.

\begin{example}
{\rm  Let $\EE$ be a nodal Weierstra\ss{} curve.
\begin{itemize}
\item
The following triple $(\widetilde\kF, \kM, \tilde{i})$
 defines an indecomposable
 vector bundle of rank 2 on
$\EE$:
the normalization $\widetilde\kF = \widetilde{\kO}\oplus \widetilde{\kO}(n), n\ne 0$, $\kM = \kk^2(s)$
and  matrices:
%my change:
%$$
%i(0) =
%\left(
%\begin{array}{cc}
%1 & 0 \\
%\hline
%0 & 1
%\end{array}
%\right), \qquad
%i(\infty) =
%\left(
%\begin{array}{cc}
%0      & 1 \\
%\hline
%\lambda & 0
%\end{array}
%\right),
%\quad \lambda \in \kk.
%$$

$$
i(0) =
\setlength\arrayrulewidth{0.3mm}
\begin{array}{
@{}c@{}     c       @{}c@{}     c       @{}c@{} l}
\cline{1-4}
%    & \hr && \hr \\
\fvr  & {\mk{1}}  &&    {\mk{0}}   & \fvr    &{~\scriptstyle 0} \\
\cline{1-4}
\fvr&  0             & &       \mk{1} & \fvr & {~\scriptstyle n}\\
%\cline{3-4}
\cline{1-4}
\end{array}
%\hspace{1cm}
~~~\hbox{and}~~
i(\infty) =
\begin{array}{
@{}c@{}     c       @{}c@{}     c       @{}c@{} l}
\cline{1-4}
%    & \hr && \hr \\
\fvr  & {\mk{0}}  &&    {\mk{1}}   & \fvr     & {~\scriptstyle 0}\\
\cline{1-4}
\fvr&  \lambda            & &       \mk{0} & \fvr & {~\scriptstyle n}  \\
%\cline{3-4}
\cline{1-4}
\end{array}
\quad \lambda \in \kk^*.
$$

\item
The triple
$\big(\widetilde{\kO}(-1), \kk^2, \tilde{i}=\framebox{1~~0}~\framebox{0~~1} \big)$
%my change:
%$(\widetilde{\kO}(-1), \kk^2, \tilde{i})$ with $\tilde{i}=i(0), i(\infty)=\framebox{1~~0}, \framebox{0~~1}$
 describes the unique torsion free sheaf that is
not locally free of degree zero, which compactifies  the Jacobian $\Pic^0(\EE)$.
%my change:
%compactifying the Jacobian $\Pic^0(\EE)$.
\end{itemize}
}
\end{example}

A Gelfand matrix problem is determined by a  certain partially ordered set together
with an  equivalence relation on it.  Such a poset with an equivalence relation is
called a \emph{bunch of chains}. Before giving
a general definition, we give an example describing the matrix problem  which
corresponds to   a nodal Weierstra\ss{} curve.

 There are two infinite
sets $E_{0}= \{E_{0}(k)| k \in {\mathbb Z} \}$ and
$E_{\infty}= \{E_{\infty}(k)| k \in {\mathbb Z} \}$ with the ordering
$\dots < E_*(-1) < E_*(0) < E_*(1) < \dots, * \in \{0,\infty\}$
and two one-point sets
$\{F_{0}\}$ and $\{F_{\infty}\}$. On the union
$$
{\mathbf E} \bigcup {\mathbf F} = (E_{0}\cup E_{\infty}) \bigcup
(F_{0} \cup F_{\infty})
$$
we introduce an equivalence relation:
$E_{0}(k)\sim E_{\infty}(k),$ where $k\in {\mathbb Z}$ and $F_{0}\sim F_{\infty}$.

%my picture
\vspace{-0.3cm}
$$
\xymatrix @ -1pc
{
&\ar@{}[r]^{\scriptstyle F_0}     &{\diamond \ar@{}[rr]}&\dhr&
 {\diamond \ar@{}[r]^{\scriptstyle F_\infty}}&&\\
                        &&{\circ\ar[d] }&\dhr& {\circ \ar[d]}&\\
\ar@{}[rr]^{\scriptstyle E_0(k)}    &&{\circ\ar[d] }&\dhr& {\circ \ar@{}[rr]^{\scriptstyle E_\infty (k)}\ar[d]}&&\\
\ar@{}[rr]^{\scriptstyle E_0 (k+1)~} &&{\circ \ar[d]}&\dhr& {\circ\ar@{}[rr]^{~\scriptstyle E_\infty (k+1)}\ar[d]}&&\\
 && && && \\}
$$

This picture contains complete
 information about the  corresponding matrix problem.
The circles denote the elements of ${\mathbf E}$, the diamonds denote
elements $F_0$ and $F_\infty,$ dotted lines connect equivalent
elements and vertical arrows describe the partial order in $E_0$ and $E_\infty.$
 The sets $E_{0}\cup F_{0}$ and  $E_{\infty}\cup F_{\infty}$ correspond to matrices $i(0)$ and
$i(\infty)$ respectively, elements $E_0(k)$ and  $E_\infty(k)$ label
 their horizontal stripes, $k \in \mathbb{Z}$.   We also say that a row from the horizontal block
$E_*(k)$ has weight $k$, where $*$ is either  $0$ or $\infty$.
The equivalence relation
$E_0(k) \sim E_\infty(k)$ means that horizontal blocks of weight $k$ have the same number
of rows and  $F_{0}\sim F_{\infty}$ tells that $i(0)$ and $i(\infty)$ have the same number of
columns. Moreover, elementary transformations inside of two conjugated blocks have to be done
simultaneously. The total order on $E_0$ and $E_\infty$ means that we can add any scalar
multiple of any row of smaller weight to any row with a bigger weight. Such transformations can
be performed in the matrices $E_0$ and $E_\infty$ independently.

\begin{definition}\label{def:cycles}
Let $\kk$ be an arbitrary field. A cycle of $n$ projective lines is a projective curve $\EE_n$
over $\kk$ with $n$ irreducible components, each  of them is isomorphic to $\mathbb{P}^1$. We
additionally assume  that all components intersect transversally with
 intersection matrix of type
$\widetilde{A}_n$ and the completion of the local rings of any
 singular point of $\EE_n$ is isomorphic to $\kk\llbracket u,v\rrbracket  /uv$.  In a similar way,
a chain of $k$ projective lines $\mathbb{I}_k$ is a configuration of projective $k$
lines with  intersection matrix of type $A_{k-1}$.
\end{definition}

\begin{remark}
{\rm
Let $\kk = \mathbb{R}$ and $\EE$ be the cubic curve $zy^2 = x^3 - zx^2$. Then $s=(0:0:1)$ is the
singular point of $\EE$ and $\widehat{\kO_{\EE,s}} = \mathbb{R}\llbracket u,v\rrbracket  /(u^2 + v^2)$.
Then $\EE$ is not a cycle  of projective lines in the sense of  Definition
\ref{def:cycles}  and the combinatorics of the indecomposable vector bundles on $\EE$
will be considered elsewhere.
}
\end{remark}

\medskip
Let
$\EE$ be either a chain or a cycle of
projective lines and $\widetilde\EE$ its
normalization.   The matrix
problem we get  is given by the following  partially ordered set.

Consider the set of pairs $\{(\mathbb{L},a) \}$, where $\mathbb{L}$ is an
irreducible component of $\widetilde\EE$ and  $a\in \mathbb{L}$ a preimage of a
singular point. To each such  pair
we attach a totally ordered set $E_{(\mathbb{L},a)} =
\{ E_{(\mathbb{L},a)}(i)| i \in {\mathbb Z} \}$,
where $ \dots < E_{(\mathbb{L},a)}(-1) < E_{(\mathbb{L},a)}(0) < E_{(\mathbb{L},a)}(1)< \dots $
and a one-point set $F_{(\mathbb{L},a)}$. On the union
$$
{\mathbf E}\bigcup {\mathbf F} = \bigcup\limits_{(\mathbb{L},a)} (E_{(\mathbb{L},a)}
\cup F_{(\mathbb{L},a)}),
$$
we introduce an equivalence relation:
\begin{enumerate}
\item $F_{(\mathbb{L}',a')}\sim F_{(\mathbb{L}'',a'')}$,
where  $a'$ and $a''$ are preimages of the same singular point $a \in \EE.$
\item $E_{(\mathbb{L},a')}(k)\sim F_{(\mathbb{L},a'')}(k)$ for $k\in {\mathbb Z}$ and
 $a', a'' \in \mathbb{L}$.
\end{enumerate}

Such a partially ordered set with an equivalence relation is
called a bunch of chains \cite{mp}. A representation of such a bunch of chains
is given by a set of matrices $M(\mathbb{L},a),$ for each element $(\mathbb{L},a).$
Every  matrix $M(\mathbb{L},a)$
 is divided into
horizontal blocks labelled by the elements of  $E_{(\mathbb{L},a)}$. Of course, all
but finitely many labels corresponds to empty blocks.
The principle of conjugation of blocks is the same as for a rational curve with
one node.

The category of representations of a bunch of chains is additive and
has two types of
indecomposable representations: \emph{bands} and \emph{strings} \cite{mp}.
Hence, the technique of representations
of bunches of chains allows to describe indecomposable torsion free sheaves on chains and
cycles of projective lines.

\medskip

Let $\EE = \EE_n$ be  a cycle of $n$
projective lines, $\{ a_1,a_2,\dots,a_n \}$ the  set of singular points
of $\EE$,
$\widetilde{\EE}\stackrel{p}\lar \EE$ the normalization of $\EE$,
$\widetilde{\EE} = \coprod\limits_{i=1}^{n} \LL_{i}$, where each $\LL_i$ is
isomorphic to a projective line  and $\{a'_{i}, a''_{i} \} = p^{-1}(a_{i})$. Assume
 that  $a'_{i},a''_{i+1} \in \LL_{i}$, where $a''_{n+1} = a''_{1}$. Fix coordinates on each projective
line  $ \LL_{i}$ in such a way, that $a'_{i} = (0: 1)$ and $a''_{i+1} = (1: 0)$.

\begin{definition}
A band $\kB({\boldsymbol  d},m, p(t))$ is  an indecomposable  vector
bundle of rank $rmk$. It is  determined  by the
following parameters:
\begin{enumerate}
\item ${\boldsymbol d}= (d_{1},d_{2},\dots,d_{n},d_{n+1},d_{n+2},\dots, d_{2n},\dots,
d_{rn-n+1}, d_{rn-n+2}, \dots,  d_{rn}) \in {\mathbb Z}^{rn}$
is a sequence of degrees on the normalized
curve $\widetilde{\EE}$. This sequence should be non-periodic, i.e. not  of the form
${\mathbf e}^s = \underbrace{ee\dots e}_{s \, times }$, where
${\boldsymbol  e} = e_1,e_2, \dots, e_{qn}$
is another sequence and  $q = \frac{r}{s}$.
\item $p(t) = t^k + a_1 t^{k-1}+\dots + a_k \in \kk[t]$ is an irreducible polynomial of
      degree $k$, $p(t) \ne t$.
\item $m \in \mathbb{Z}^+$ is a positive integer.
\end{enumerate}

\noindent
In particular, one can recover  from the sequence $\boldsymbol{d}$
the pull-back   of $\kB({\boldsymbol  d},m, p(t))$
on the $l$-th irreducible component of $\widetilde\EE$: it  is
$$p_l^*(\kB({\boldsymbol  d},m, p(t))) \cong \bigoplus\limits_{i=1}^{r} \kO_{\LL_{l}}(d_{l+in})^{mk}.$$

\medskip
\noindent
A string $\kS(\boldsymbol{d}, f)$ is a torsion free sheaf which
depends only on two discrete parameters $f\in \{1,2,\dots,n\}$ and
$\boldsymbol{d} = (d_1, d_2,\dots, d_t)$, $t>1$.
\end{definition}

Now we are going to explain the way of construction of gluing matrices of triples
corresponding to bands $\kB(\boldsymbol{d}, m, p(t))$ and strings
$\kS(\boldsymbol{e},f)$.

\begin{algorithm}\label{alg:canf}
{\rm
\medskip
\textbf{Bands}. Let
$\boldsymbol{d}= (d_1, d_2,\dots, d_{rn}) \in \mathbb{Z}^{rn}$ be a non-periodic sequence,
$m \in \mathbb{Z}^+$ and $p(t) \in \kk[t]$ an irreducible polynomial of degree $k$.
We have $2n$ matrices $M(\LL_{i},a'_{i})$
and $M(\LL_{i},a''_{i+1}), i=1,\dots,n$
occurring in the triple, corresponding to $\kB({\boldsymbol  d}, m, p(t))$. Each of them has
size $mrk\times mrk$. Divide these matrices into $mk\times mk$ square blocks.
Consider the  sequences ${\boldsymbol  d}(i)=d_{i}d_{i+n}\dots d_{i+(r-1)n}$ and label
the horizontal strips of $M(\LL_{i}, a'_{i})$ and $M(\LL_{i},a''_{i+1})$ by
integers occurring  in each  ${\boldsymbol d}(i)$. If some
integer $d$ appears  $l$ times in ${\boldsymbol  d}(i)$ then the horizontal strip
corresponding  to the label $d$ consists of $l$ substrips having   $mk$ rows each.
Recall now an algorithm of writing the components of the matrix $\tilde{i}$ in a
normal form \cite{mp}:

\begin{enumerate}
\item Start with the  sequence
$
(\LL_{1},a'_{1})\stackrel{1}\lar (\LL_{1},a''_{2})\stackrel{1}\lar (\LL_{2},a'_{2})\stackrel{1}\lar
(\LL_{2},a''_{3})\stackrel{1}\lar \dots \stackrel{1}\lar
(\LL_{n},a'_{1}) \stackrel{2}\lar (L_{1},a''_{1})\stackrel{2}\lar \dots
 \stackrel{r}\lar (\LL_{n},a'_{n})
\stackrel{r}\lar (\LL_{n},a''_{1}) \stackrel{1}\lar.
$
It is convenient to imagine  this sequence as a cyclic word   broken at the place
$(\LL_{1},a'_{1}).$

\item Unroll the sequence ${\boldsymbol d}$. This means that we write over each
$(\LL_{i},a)$ the corresponding term of the subsequence ${\boldsymbol  d(i)}$
together with the  number of its previous  occurrences  in ${\boldsymbol  d(i)}$
including the current one:
\begin{align*}
(\LL_{1},a'_{1})^{(d_{1},1)}\stackrel{1}\lar (\LL_{2},a'_{2})^{(d_{2},1)}
\stackrel{1}\lar
(\LL_{2},a''_{3})^{(d_{2},1)} \stackrel{1}\lar \\
\dots\stackrel{r}\lar (\LL_{n}, a''_{1})^{d_{rn},*}
\stackrel{r}\lar
(\LL_{n}, a'_{1})^{(d_{rn},1)} \stackrel{1}\lar.
\end{align*}

\item Now we can  fill the entries of the matrices $M(\LL,a)$. Consider each arrow \\
$
(\LL, a)^{(d,i)} \stackrel{l}\lar.
$
Then insert  the matrix $I_{mk}$ in the block $((d,i),l)$ of the matrix
$M(\LL,a)$, which is defined as the  intersection of the
$i$-th substrip of the horizontal strip labeled by  $d$
and the $l$-th vertical strip.
\item Put at the $((d_{rn},1),r)$-th place of $M(\LL_{n},a''_{1})$ the
 Frobenius block $J_m(p(t))$.
\end{enumerate}

\noindent
\textbf{Strings}.
Let $\boldsymbol{e} = (e_1,e_2,\dots,e_s) \in \mathbb{Z}^s$ and
$f \in \{1,2,\dots,n\}$.
 The algorithm to write
the matrices for  the torsion free sheaf $\kS(\boldsymbol{d},f)$ is essentially the
same as for bands.
The parameter $f$ denotes the number of the component $\mathbb{L}_f$ of $\widetilde\EE$
which
the sequence
$
(\LL_{f},a'_{f}) \lar (\LL_{f},a''_{f+1}) \lar
\dots \lar
(\LL_{f+s},a'_{f+s}) \lar (\LL_{f+s},a''_{f+s+1})
$
starts with.
Then we unroll $\boldsymbol{e}$  using the same algorithm as for bands.
The only difference is that we insert instead of $I_{mk}$ the unit  $(1\times 1)$--matrix.
Note, that some of the matrices $M(\mathbb{L}_i, a_i')$ and $M(\mathbb{L}_i, a_i'')$ can be
 non-square, but they are automatically of full row rank.
}
\end{algorithm}

\begin{example}
{\rm
Let  $\EE = \EE_2$ be a cycle of two projective lines,
$\boldsymbol{d} = (0,1,1,3,1,-2)$ and $p(t) \in \kk\llbracket t\rrbracket  $ an irreducible polynomial of degree
$k$.
 Then  ${\mathcal B}({\boldsymbol d},m, p(t))$
is a vector bundle of degree $3m$
with the normalization
$$\bigl({\mathcal O}_{\LL_{1}}^{mk} \oplus {\mathcal O}_{\LL_{1}}(1)^{2mk}\bigr) \oplus
\bigl({\mathcal O}_{\LL_{2}}(-2)^{mk} \oplus {\mathcal O}_{\LL_2}(1)^{mk} \oplus
{\mathcal O}_{\LL_2}(3)^{mk}\bigr)$$
and  gluing matrices

%my figure:
\begin{align*}
\setlength\arrayrulewidth{0.1mm}
\setlength\doublerulesep{0.1mm}\doublerulesepcolor{black}
M(\LL_1, a_1^{\prime})=
\begin{array}{
@{}c@{}     c       @{}c@{}     c       @{}c@{}c       @{}c@{} }
\hline\hline
%    & \hr && \hr \\
\fvr & {\mk{\scriptstyle I_{mk}}}  &\vr&    \mk{}   & \vr  &{\mk{}}   & \fvr   \\
\hline\hline
\fvr & {\mk{}}            &\vrule&        \mk{\scriptstyle I_{mk}} & \vr  &{\mk{}}   & \fvr \\
\cline{1-6}
\fvr & {\mk{}}            &\vrule&        \mk{} & \vr  &{\mk{\scriptstyle I_{mk}}}   & \fvr \\
\hline\hline
\end{array}
\begin{array}{r}
{\scriptstyle 0}\\
 \\
{^1}\\
\end{array}
~
\begin{array}{
@{}c@{}     c       @{}c@{}     c       @{}c@{}c       @{}c@{} }
\hline\hline
%    & \hr && \hr \\
\fvr & \mk{\scriptstyle I_{mk}}  &\vr&    \mk{}   & \vr  &{\mk{}}   & \fvr   \\
\hline\hline
\fvr & {\mk{}}            &\vrule&        \mk{\scriptstyle I_{mk}} & \vr  &{\mk{}}   & \fvr \\
\cline{1-6}
\fvr & {\mk{}}            &\vrule&        \mk{} & \vr  &{\mk{\scriptstyle I_{mk}}}   & \fvr \\
\hline\hline
\end{array}
=
M(\LL_1, a_2^{\prime\prime})
\end{align*}
%\\
\begin{align*}
\setlength\arrayrulewidth{0.1mm}
\setlength\doublerulesep{0.1mm}\doublerulesepcolor{black}
M(\LL_2, a_1^{\prime\prime})=
\begin{array}{
@{}c@{}     c       @{}c@{}     c       @{}c@{}c       @{}c@{} }
\hline\hline
%    & \hr && \hr \\
\fvr & {\mk{\scriptstyle J_{m}}}  &\vr&    \mk{}   & \vr  &{\mk{}}   & \fvr     \\
\hline\hline
\fvr & {\mk{}}            &\vrule&        \mk{\scriptstyle I_{mk}} & \vr  &{\mk{}}   & \fvr \\
%\cline{1-6}
\hline\hline
\fvr & {\mk{}}            &\vrule&        \mk{} & \vr  &{\mk{\scriptstyle I_{mk}}}   & \fvr   \\
\hline\hline
\end{array}
\begin{array}{r}
{\rule{0pt}{12pt}\scriptstyle -2}\\
{\rule{0pt}{12pt}\scriptstyle 1}\\
\rule{0pt}{12pt}\scriptstyle 3\\
\end{array}
~
\begin{array}{@{}c@{}     c       @{}c@{}     c       @{}c@{}c       @{}c@{} }
\hline\hline
%    & \hr && \hr \\
\fvr & {\mk{}}  &\vr&    \mk{}   & \vr  &{\mk{\scriptstyle I_{mk}}}   & \fvr   \\
\hline\hline
\fvr & {\mk{\scriptstyle I_{mk}}}            &\vrule&        \mk{} & \vr  &{\mk{}}   & \fvr \\
%\cline{1-6}
\hline\hline
\fvr & {\mk{}}            &\vrule&        \mk{\scriptstyle I_{mk}} & \vr  &{\mk{}}   & \fvr  \\
\hline\hline
\end{array}
=
M(\LL_2, a_2^{\prime}),
\end{align*}
where $J_m$ is the Frobenius block corresponding to the $\kk[t]$-module  $\kk[t]/p(t)^m.$

\noindent
The corresponding unrolled sequence looks as follows:

\noindent
$%\begin{eqnarray*}
(\LL_{1},a'_{1})^{(0,1)} \stackrel{1}\longrightarrow
(\LL_{1},a''_{2})^{(0,1)}\stackrel{1}\longrightarrow
(\LL_{2},a'_{2})^{(1,1)}\stackrel{1}\longrightarrow
(\LL_{2},a''_{1})^{(1,1)}
\stackrel{2}\longrightarrow
\\
(\LL_{1},a'_{1})^{(1,1)} \stackrel{2}\longrightarrow
(\LL_{1},a''_{2})^{(1,1)}\stackrel{2}\longrightarrow
(\LL_{2},a'_{2})^{(3,1)}\stackrel{2}\longrightarrow
(\LL_{2},a''_{1})^{(3,1)}
\stackrel{3}\longrightarrow
\\
(\LL_{1},a'_{1})^{(1,2)}\stackrel{3}\longrightarrow
(\LL_{1},a''_{2})^{(1,2)}\stackrel{3}\longrightarrow
(\LL_{2},a'_{2})^{(-2,1)}\stackrel{3}\longrightarrow
(\LL_{2},a''_{1})^{(-2,1)}\stackrel{1}\lar.
$
%\end{eqnarray*}

\medskip
\medskip
\noindent
Let $f=2$ and $\boldsymbol{d} = (-1,0,1,-1,1)$.
Then the corresponding torsion free sheaf $\kS(\boldsymbol{d},f)$ has normalization
$$
\widetilde\kF = \bigl(\kO_{\LL_1}(-1) \oplus \kO_{\LL_1}\bigr) \oplus
\bigl(\kO_{\LL_2}(-1) \oplus \kO_{\LL_2}(1)^2\bigr)
$$
and gluing matrices

%my picture
\begin{align*}
\setlength\arrayrulewidth{1pt}
\setlength\doublerulesep{0.1pt}\doublerulesepcolor{black}
M(\LL_1, a_1^{\prime})=
\begin{array}{
@{}c@{}     c       @{}c@{}     c       @{}c@{}c       @{}c@{} }
\cline{1-6}
%    & \hr && \hr \\
\fvr & {\mk{}}  &\vr&    \mk{1}   & \vr  &{\mk{}}   & \fvr   \\
\cline{1-6}
\fvr & {\mk{1}}            &\vrule&        \mk{} & \vr  &{\mk{}}   & \fvr \\
\cline{1-6}
\end{array}
\begin{array}{r}
{\scriptstyle -1}\\
{\scriptstyle 0}
\end{array}
~
\begin{array}{@{}c@{}     c       @{}c@{}     c       @{}c@{}c       @{}c@{} }
\cline{1-6}
%    & \hr && \hr \\
\fvr & {\mk{}}  &\vr&    \mk{}   & \vr  &{\mk{1}}   & \fvr   \\
\cline{1-6}
\fvr & {\mk{}}            &\vrule&        \mk{1} & \vr  &{\mk{}}   & \fvr \\
\cline{1-6}
\end{array}
=
M(\LL_1, a_2^{\prime\prime})
\end{align*}
%\\
\begin{align*}
\setlength\arrayrulewidth{0.1mm}
\setlength\doublerulesep{0.1mm}\doublerulesepcolor{black}
M(\LL_2, a_1^{\prime\prime})=
\begin{array}{
@{}c@{}     c       @{}c@{}     c       @{}c@{}c       @{}c@{} }
\hline\hline
%    & \hr && \hr \\
\fvr & {\mk{1}}  &\vr&    \mk{}   & \vr  &{\mk{}}   & \fvr     \\
\hline\hline
\fvr & {\mk{}}            &\vrule&        \mk{1} & \vr  &{\mk{}}   & \fvr \\
\cline{1-6}
\fvr & {\mk{}}            &\vrule&        \mk{} & \vr  &{\mk{1}}   & \fvr   \\
\hline\hline
\end{array}
\begin{array}{r}
{\scriptstyle -1}\\
{\rule{0pt}{18pt} \scriptstyle 1}\\
\,
\end{array}
~
\begin{array}{@{}c@{}     c       @{}c@{}     c       @{}c@{}c       @{}c@{} }
\hline\hline
%    & \hr && \hr \\
\fvr & {\mk{1}}  &\vr&    \mk{}   & \vr  &{\mk{}}   & \fvr   \\
\hline\hline
\fvr & {\mk{}}            &\vrule&        \mk{1} & \vr  &{\mk{}}   & \fvr \\
\cline{1-6}
\fvr & {\mk{}}            &\vrule&        \mk{} & \vr  &{\mk{1}}   & \fvr  \\
\hline\hline
\end{array}
=
M(\LL_2, a_2^{\prime})
\end{align*}

\noindent
The corresponding unrolled sequence is

\noindent
$
(\LL_{2},a'_{2})^{(-1,1)}\stackrel{1}\longrightarrow
(\LL_{2},a''_{1})^{(-1,1)}\stackrel{1}\longrightarrow
(\LL_{1},a'_{1})^{(0,1)}\stackrel{1}\longrightarrow
(\LL_{1},a''_{2})^{(0,1)}
\stackrel{1}\longrightarrow\\
(\LL_{2},a'_{2})^{(1,1)}\stackrel{2}\longrightarrow
(\LL_{2},a''_{1})^{(1,1)}\stackrel{2}\longrightarrow
(\LL_{1},a'_{1})^{(-1,1)}\stackrel{2}\longrightarrow
(\LL_{1},a''_{2})^{(-1,1)}
\stackrel{2}\longrightarrow\\
(\LL_{2},a'_{2})^{(1,2)}\stackrel{3}\longrightarrow
(\LL_{2},a''_{1})^{(1,2)}\stackrel{3}\longrightarrow.
$
}
\end{example}

\noindent
Summing everything up, we get  the following theorem.

\begin{theorem}[\cite{DrozdGreuel}]
Let $\EE = \EE_n$ be a cycle of $n$ projective lines over a filed $\kk$.
Then
\begin{itemize}
\item
any indecomposable vector bundle on $\EE$
is isomorphic to some $\kB(\boldsymbol{d}, m, p(t))$, where  $\boldsymbol{d} =
(d_1,d_2,\dots, d_{rn}) \in \mathbb{Z}^{rn}$ is a non-periodic sequence,
 $m \in \mathbb{Z}^+$ and
$p(t) = t^k + a_1 t^{k-1} + \dots + a_k \in \kk[t]$
is an irreducible polynomial, $p(t) \ne t$.
\item
Any  torsion free but not locally
free coherent sheaf is isomorphic to some $\kS(\boldsymbol{d}, t)$, where
$t \in \{1,2,\dots,n\}$ and $\boldsymbol{d} \in \mathbb{Z}^{rn}$.
\end{itemize}
The only isomorphisms between
indecomposable vector bundles are generated by
\begin{itemize}
\item $\kB(\boldsymbol{d}, m, p(t)) \cong \kB(\boldsymbol{d}^{\circ}, m, q(t)),$
$\boldsymbol{d}^\circ = d_{rn},d_{rn-1},\dots,d_1$ and $q(t) = \frac{t^k}{a_k} p(1/t)$
\item $\kB(\boldsymbol{d}, m, p(t)) \cong \kB(\boldsymbol{d'}, m, p(t))$, with
$$\boldsymbol{d'} = d_{n+1},d_{n+2},\dots, d_{2n},d_{2n+1},\dots, d_1,d_2,\dots,d_n.$$
\end{itemize}
The only isomorphisms between strings are
$
\kS(\boldsymbol{e},f) \cong \kS(\boldsymbol{e}^\circ,f^\circ),
$
where $\boldsymbol{e}^\circ$ is the opposite sequence
$\boldsymbol{e}^\circ = e_s, e_{s-1}, \dots,e_1$.
If $s = nk + s'$ with  $0 \le s' < n$,  then
$f^{\circ} =  s' + f$ taken modulo  $n$.
\end{theorem}

\begin{remark}
If $\EE = \EE_1$ is a Weierstra\ss{} nodal curve, then there is no choice for
the parameter $f$ in the definition of a string and one simply uses the notation
$\kS(\boldsymbol{d})$. If the field $\kk$ is algebraically closed, we write
$\kB(\boldsymbol{d},m,\lambda)$ instead of $\kB(\boldsymbol{d},m, t-\lambda)$.
\end{remark}

As a direct corollary of the combinatorics of strings and bands
we obtain  the following theorem

\begin{theorem}[\cite{DrozdGreuel}]
Let $\XX = \mathbb{I}_n$ be a chain of $n$ projective lines, then any vector bundle
$\kE$ on $\XX$ splits into a direct sum of lines bundles.  Moreover,
$\Pic(\mathbb{I}_n) \cong
\mathbb{Z}^n$ and  a line bundle is determined by its restrictions on each irreducible component.
A description of torsion free sheaves is similar: any indecomposable torsion free sheaf is
isomorphic to  the direct image of a line bundle on a subchain of projective lines.
\end{theorem}

\subsection{Properties of torsion free sheaves
on cycles of projective lines}\label{subsec:prop}
~~

\noindent
Throughout this subsection,  let $\kk$ be an algebraically closed field and
 $\EE = \EE_n$ a cycle of $n$ projective lines over $\kk$.
As we have seen in the previous subsection, indecomposable  vector bundles
on $\EE$ are bands $\kB(\boldsymbol{d}, m,\lambda)$ and indecomposable torsion
free but  not locally sheaves are strings $\kS(\boldsymbol{d}, f)$. They were described
in terms of a certain problem of linear algebra. However, in the case of an algebraically
closed field
there is   a geometric way to present the classification of  indecomposable torsion free
sheaves on $\EE$ without appealing to the formalism of bunches of chains.
This description, the proof of which we give here for the first time,
 is completely parallel to Oda's one for vector bundles
on elliptic curves \cite{Oda}.

We start with a lemma describing unipotent vector bundles on $\EE$.

\begin{lemma}
For any integer $m \ge 1$ there exists a unique indecomposable  vector bundle
$\kF_m$  on $\EE_n$ appearing in the exact sequence
$$
0 \lar \kF_{m-1} \lar \kF_m \lar \kO \lar 0, \qquad \kF_1 = \kO.
$$
In our notation we have $\kF_m \cong \kB(\boldsymbol{0},m,1)$, where
$\boldsymbol{0} = \underbrace{(0,0,\dots,0)}_{n \, times}$.
\end{lemma}

\noindent
\emph{Sketch of the proof}.
Since the dualizing sheaf $\omega_\EE \cong \kO$ is trivial, we have
$$\Ext^1(\kO, \kO) = \kk$$ and there is a unique non-split extension
$$
0 \lar \kO \lar \kF_2 \lar \kO \lar 0.
$$
Then using the same arguments as in \cite{Atiyah}, we can inductively
construct
indecomposable
vector bundles $\kF_m$, $m \ge 1$ such that $H^0(\kF_m) \cong H^1(\kF_m) = \kk$
together with  exact sequences
$$
0 \lar \kF_{m-1} \lar \kF_m \lar \kO \lar 0, \qquad \kF_1 = \kO.
$$
On the other hand, $\kB(\boldsymbol{0},m,1)$ is the unique indecomposable vector bundle
on $\EE$ of rank $m$ and normalization $\kO_{\widetilde\EE}^m$
with  a non-zero section.
Hence $\kF_m \cong \kB(\boldsymbol{0},m,1)$.

%\medskip
\medskip
\noindent
The proof of the following proposition is straightforward:

\begin{proposition}
Let $\Psi: \VB(\EE) \lar \mathsf{T}_\EE$ be the functor establishing an
equivalence between the category of vector
bundles on $\EE$ and the category of triples.
Then $\Psi$ preserves tensor products:
$
\Psi(\kE\otimes \kF) \cong \Psi(\kE)\otimes \Psi(\kF),
$
where $(\widetilde\kE, \kM, \tilde{i}) \otimes (\widetilde\kF, \kN, \widetilde{j}) =
(\widetilde\kE \otimes_{\widetilde\kO}  \widetilde\kF, \kM\otimes_\kA \kN, \tilde{i}\otimes \widetilde{j}).$
In particular,
\begin{itemize}
\item We have an isomorphism
$
\kB((d), m,\lambda) \cong \kB((d), 1,\lambda)\otimes \kF_m.
$
\item There is the following rule for a decomposition of
the tensor product of two unipotent vector bundles:
$$
\kF_f \otimes \kF_g \cong \bigoplus\limits_{i} \kF_{h_i},
$$
where integers $h_i$ are the same as in the  decomposition
$$\kk[t]/t^f \otimes_{\kk} \kk[t]/t^g \cong
\bigoplus\limits_{i \in \mathbb{Z}} \kk[t]/t^{h_i}$$
 in the category of $\kk[t]$--modules.
\item In particular, if $\kk$ is of characteristics zero, we have
$$
\kF_f \otimes \kF_g  \cong
\bigoplus\limits_{j=1}^{g} \kF_{f-g-1+2j}.
$$
\end{itemize}
\end{proposition}

Now we formulate a geometric description of indecomposable torsion free sheaves
on a cycle of projective lines in the case of an algebraically closed field.

\begin{theorem}\label{thm:etale}
Let $\EE = \EE_n$ be a cycle of $n$ projective lines and $\mathbb{I}_k$  be a  chain of
$k$ projective lines, $\kE$ an indecomposable   torsion free sheaf on $\EE_n$.
\begin{enumerate}
\item If $\kE$ is  locally free, then there is an \'etale covering
$\pi_r: \EE_{nr}\lar \EE_n$, a line bundle $\kL\in \Pic(\EE_{nr})$ and
a natural number $m\in {\mathbb N}$ such that
$$
 \kE \cong \pi_{r*}(\kL \otimes \kF_m).
$$

%\vspace{-1cm}
%insert picture
\begin{figure}[ht]
\hspace{-1cm} %\vspace{-1cm}
\psfrag{A}{$\EE_4$}
\psfrag{B}{$\EE_2$}
\psfrag{p}{$\pi_2$}
\includegraphics[height=2cm,width= 5cm,angle=0]{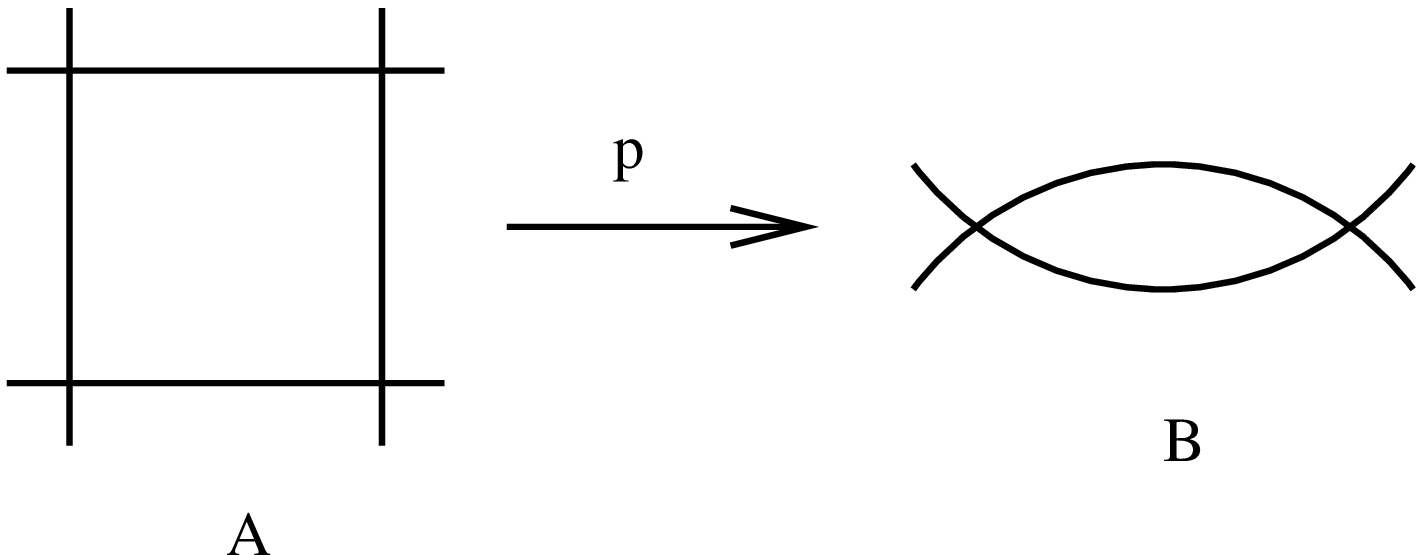}
\end{figure}

\noindent
Moreover, if $char(\kk) = 0$, then integers $r$ and $m$ are uniquely determined and
the line bundle $\kL$ is unique up to the action of ${\mathrm{Aut}}(\EE_{nr}/\EE_n)$. 
Other way arround, for 
given $r, m \in \mathbb{Z}^+$ and $\kL \in \Pic(\EE_{nr})$ 
the vector bundle $\pi_{r*}(\kL \otimes \kF_m)$ is indecomposable 
if and only if $\kL$ does not belong to the image of the map 
$\pi^*_t: \Pic(\EE_{nr/t}) \lar  \Pic(\EE_{nr})$ for any proper 
divisor $t$  of 
$r$. 

\item If $\kE$ is not locally free then there exists a map
$p_k: \mathbb{I}_{k} \lar \EE_n$ and a uniquely determined 
line bundle $\kL \in \Pic(\mathbb{I}_k)$ such that
$\kE \cong p_{k*}(\kL)$.
%insert picture
\begin{figure}[ht]
\hspace{-1cm}
\psfrag{A}{$\mathbb{I}_{4}$}
\psfrag{B}{$\EE_1$}
\psfrag{p}{$p_4$}
\includegraphics[height=2cm,width=4cm,angle=0]{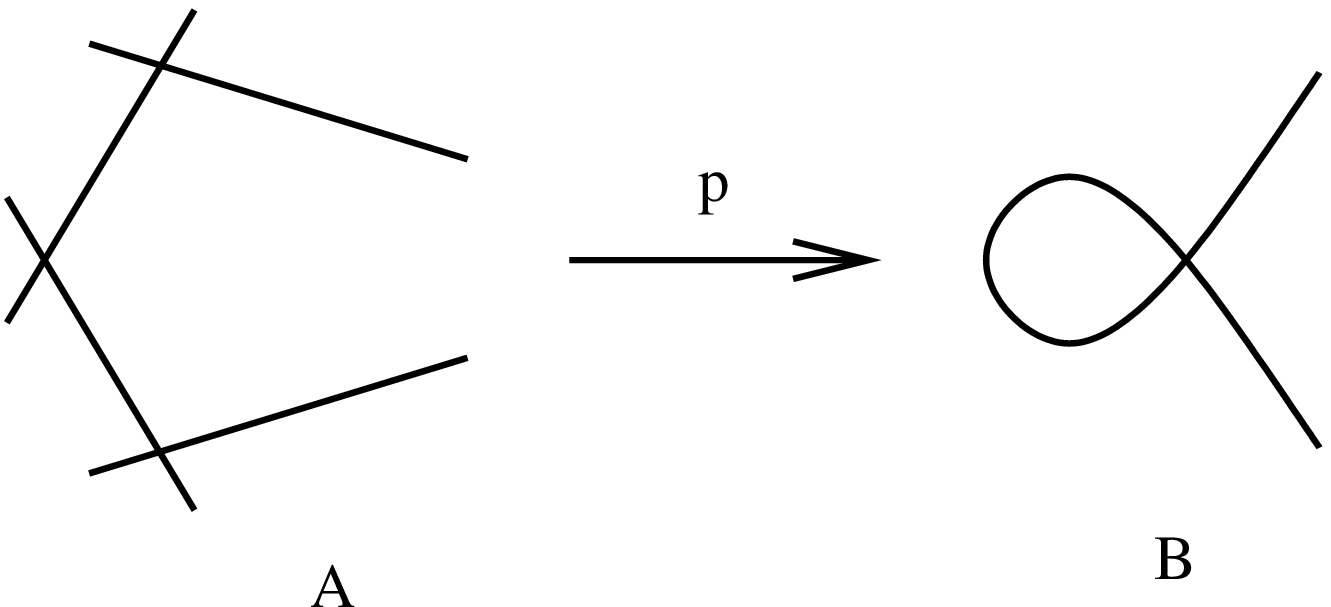}
\end{figure}

\noindent
Other way arround, the torsion free sheaf $p_{k*}(\kL)$ is indecomposable for any line
bundle $\kL \in \Pic(\mathbb{I}_k)$. 
\end{enumerate}
\end{theorem}

\noindent
\emph{Proof}. Let $\EE$ be a cycle of projective lines and
$\pi_\EE: \EE' \lar \EE$ an \'etale covering of degree $r$, $\kF$ a
torsion free sheaf on $\EE'$. In the notation of Remark \ref{rem:triples}
we have the commutative  diagram
$$
  \xymatrix{
   &  \SS' \ar[rrr]^{\pi_{\SS}} \ar@{-}[d]^{i'}  & &  & \SS  \ar[dd]^i  \\
   \widetilde{\SS'} \ar[rrr]^{\pi_{\widetilde{\SS}}} \ar[dd]_{\widetilde{i'}}
\ar[ru]^{\widetilde{p'}}  & \ar [d]  &   &  \widetilde\SS \ar[dd]\ar@{}[d]^{\tilde{i}} \ar[ru]_{\tilde{p}}  &  \\
   &  \EE' \ar@{-}[rr]_{\pi_{\EE}}  &  &\ar[r] &  \EE   \\
   \widetilde{\EE}' \ar[rrr]_{\pi_{\widetilde\EE}} \ar[ru]_{p'}  &  &  &
\widetilde\EE \ar[ru]_p  &
  }
$$
in which all squares are pull-back diagrams. In order to prove the theorem we have to compute
the triple describing the torsion free sheaf $\pi_*(\kF)$.  Note that  each map
$\mathbb{I} \lar \EE$
from a chain
of projective lines to a cycle of projective lines factors through an \'etale covering
$\mathbb{E}' \lar \EE$.  So, in order to prove the second part of the theorem about the
characterization of strings we may  consider an \'etale covering of $\EE$ as well.

\medskip

Note the  following simple fact about
pull-back diagrams:

\begin{lemma}
Let
$$
\xymatrix{
 \YY'  \ar[r]^{g'} \ar[d]_{f'}  & \YY \ar[d]^{f}\\
 \XX' \ar[r]^{g}  &  \XX
}
$$
be a  pull-back diagram, where all maps $f,g,f',g'$ are affine. Then for any coherent sheaf
$\kF$ on $\YY$ it holds
$
g^* f_*\kF \cong f'_* g'^* \kF.
$
\end{lemma}

\noindent
The morphism
$
p^*\pi_{\EE*}(\kF)/\mathsf{tor}(p^*\pi_{\widetilde{E}*}(\kF))
\lar \pi_{\widetilde{\EE}*}\bigl(p'^*(\kF)/\mathsf{tor}(p'^*(\kF))\bigr)
$
is an isomorphism. Indeed, we have a surjection
$$
p^*\pi_{\EE*}(\kF) \stackrel{\cong}\lar \pi_{\widetilde{\EE}*} p'^*\kF \lar
\pi_{\widetilde{\EE}*}(p'^*(\kF)/\mathsf{tor}(p'^*(\kF)))
$$
which induces a surjective map $$p^*\pi_{\EE*}(\kF)/\mathsf{tor}(p^*\pi_{\widetilde{\EE}*}(\kF))
\lar \pi_{\widetilde{\EE}*}\bigl(p'^*(\kF)/\mathsf{tor}(p'^*(\kF))\bigr)$$
of torsion free sheaves.
Since both
sheaves have the
same rank on each irreducible component of $\widetilde{\EE}$,
we conclude that this map is also injective and
therefore an isomorphism.

\medskip

We need one more simple statement about \'etale coverings.

\begin{lemma}
Let $\pi: \mathbb{Y} \lar \XX$ be an \'etale map of reduced schemes
and $\kF$ a coherent sheaf on $\mathbb{Y}$.
Then there is a canonical isomorphism 
$$\pi_*(\kF/\mathsf{tor}(\kF))
\lar \pi_*(\kF)/\mathsf{tor}(\pi_*(\kF)).$$
\end{lemma}

\noindent
\emph{Proof}. The canonical map $\kF \to \kF/\mathsf{tor}(\kF)$ induces
the morphism $\pi_*(\kF) \lar \pi_*(\kF/\mathsf{tor}(\kF))$. Since $\pi$ is \'etale,
the sheaf $\pi_*(\kF/\mathsf{tor}(\kF))$ is torsion free and the induced morphism
$$\pi_*(\kF/\mathsf{tor}(\kF))  \lar \pi_*(\kF)/\mathsf{tor}(\pi_*(\kF))$$ is an
isomorphism since it is an isomorphism on the stalks.

\medskip
\noindent
Let $\widetilde{\varepsilon}: \tilde{p}'^*i'^*(\kF) \lar
\widetilde{i'}^*(p'^*(\kF)/\mathsf{tor}(p'^*(\kF)))$
be the gluing map describing the torsion free sheaf $\kF$ in the corresponding triple.
From the commutativity of the  diagram
$$
\xymatrix{
 \tilde{p}^*i^*\pi_{\EE*}(\kF) \ar[rrrr]  \ar[dd]& & & &
\tilde{i}^*\bigl(p^* \pi_{\EE*}(\kF)/\mathsf{tor}(p^* \pi_{\EE*}(\kF))\bigr) \ar[d] \\
& & & &
\tilde{i}^*\bigl(\pi_{\widetilde{\EE}*}p'^*(\kF)/\mathsf{tor}(\pi_{\widetilde{\EE}*}p'^*(\kF))\bigr)
\ar[d] \\
\tilde{p}^*\pi_{\widetilde{\SS}*}i_n^*(\kF) \ar[d]& & & &
\tilde{i}^*\bigl(\pi_{\widetilde{\EE}*}(p'^*(\kF)/\mathsf{tor}(p'^*(\kF))\bigr) \ar[d] \\
\pi_{\widetilde{\SS}*}\widetilde{p'}^*i'^*(\kF)
\ar[rrrr]^{\pi_{\widetilde{\SS}*}(\widetilde{\varepsilon})} & & & &
\pi_{\widetilde{\SS}*}\widetilde{i'}^*\bigl(p'^*(\kF)/\mathsf{tor}(p'^*(\kF))\bigr)
}
$$
we conclude that the direct image sheaf $\pi_{\EE*}\kF$ is described by the gluing matrices
$\pi_{\widetilde{\mathbb{S}}*}(\widetilde{\varepsilon})$, which are exactly
the matrices constructed  in the Algorithm \ref{alg:canf}.
This completes the proof.

\medskip
\medskip

\begin{remark}
{\rm For a given cycle of projective lines $\EE$ over an arbitrary field $\kk$
there is always an \'etale covering  $\pi: \EE' \lar \EE$ of a given degree $r$.
For example, let $\EE = \EE_1$ be a rational curve with one node, $\XX_i = \PP^1$ ($i=1,2$),
$f_i: \Spec(\kk\times\kk) \lar \XX_i$ be two closed embeddings
 with the image $0$ and $\infty$
and $g_i: \XX_i \lar \EE$  two
normalization maps mapping the points $0$ and $\infty$ on $\XX_i$
 to the singular point of $\EE$.
  Then the push-out  of  $\XX_1$ and $\XX_2$ over $\Spec(\kk\times\kk)$
(in the category of all schemes and affine maps)
  is
a cycle of two projective lines $\EE_2$ and the induced map $g: \EE_2 \lar \EE_1$ is an \'etale
covering of degree two. The general case can be considered in a similar  way. Note that
this is quite different to the case of elliptic curves, where the existence of an \'etale
covering of a given degree strongly depends on the arithmetics of the curve.
}
\end{remark}

Similarly to the proof of Theorem \ref{thm:etale}  we have the following proposition.

\begin{proposition}
Let  $\pi_r: \EE_{nr} \lar \EE_n$
be an \'etale covering of degree $r$, $\widetilde\EE_n = \coprod\limits_{i=1}^n \LL_i$ and
$\widetilde\EE_{nr} = \coprod\limits_{j=1}^{nr} \mathbb{L}_{j}'$ be the normalizations of $\EE_n$
and
$\EE_{nr}$. Let
$\{a_1, a_2, \dots, a_n\}$ be the set
of singular points of $\EE_n$ and $\{b_1, b_2, \dots, b_n, b_{n+1}, \\\dots, b_{nr}\}$
the singular points  of $\EE_{nr}$ and  $\pi_{r}^{-1}(a_i) =
\{b_i, b_{i+n}, \dots, b_{i+ (r-1)}\}$.
Assume $\kE$ is a vector  bundle of rank $l$, given by the triple $(\widetilde\kE, \kA^l,
\tilde{i})$, where
$\widetilde\kE \cong \widetilde\kE_1 \oplus \widetilde\kE_2 \oplus \dots \oplus \widetilde\kE_n$
and $\tilde{i}$ is given by matrices
$M(\LL_1, a'_1),~ M(\LL_1, a''_2), ~\dots,\\ M(\LL_n, a'_1)$. Then the pull-back
$\pi_r^*(\kE)$ corresponds to the triple $(\widetilde{\kE}', \kA'^l, \tilde{i}'),$ where
$\widetilde{\kE}'|_{\LL'_{i+nj}} \cong \widetilde{\kE}_i$, $0 \le j \le r-1$
 and
$\widetilde{i'}$ is given by matrices $M(\LL', b') = M(\LL, a')$ and
$M(\LL', b'')  = M(\LL, a'')$ if $\pi_r(\LL') = \LL$ and  $\pi_r(b) = a$.
\end{proposition}

From this proposition follows the following corollary:

\begin{corollary}
Let Let $\EE = \EE_n$ be a cycle of $n$ projective lines and $\pi_r: \EE_{rd} \lar \EE_n$
an \'etale covering of degree $r$. Then
\begin{enumerate}
\item $\pi_r^*\kB(\boldsymbol{d}, 1, \lambda) \cong
\kB(\boldsymbol{d}^r, 1, \lambda^r)$.
\item If $char(\kk) = 0$, then $\pi_r^*(\kF_m) \cong \kF_m$. In particular,
we have an isomorphism
$$
\kB(\boldsymbol{d}, m, \lambda) \cong \kB(\boldsymbol{d}, 1, \lambda) \otimes \kF_m.
$$
\end{enumerate}
\end{corollary}

\noindent
\emph{Proof}. The proof of the first part is straightforward.
To prove the second, note that $\pi_r^*(\kF_m)$ is given by a triple, isomorphic to
$(\widetilde\kO^m, \kA^m, \tilde{i})$, where
$\tilde{i}$ is given by matrices
$M(\LL_1, a'_1) = I_m$, $M(\LL_1, a''_2) = I_m,\, \dots\,,   M(\LL_n, a'_n) = I_m,\,
M(\LL_n, a'_1) = J_m(1)^r$, where
$J_m(1)$ is the Jordan $(m\times m)$--block with the eigenvalue $1$.
If $char(\kk) = 0$ then $J_m(1)^r \sim J_m(1)$ and we get the claim. Note, that in the
case $char(\kk) = p$ we have  $J_p(1)^p = I_p$ that implies $\pi_p^*(\kF_p) \cong \kO^p$.
To complete the proof of the second claim note, that
$$
\kB(\boldsymbol{d}, m,\lambda) \cong \pi_{r*}(\kL(\boldsymbol{d}, \lambda) \otimes \kF_m)
\cong\pi_{r*}(\kL(\boldsymbol{d}, \lambda) \otimes \pi_r^*\kF_m) \cong
$$
$$
 \pi_{r*}(\kL(\boldsymbol{d}, \lambda)) \otimes \kF_m \cong \kB(\boldsymbol{d}, 1, \lambda)
\otimes \kF_m.
$$
\medskip

\begin{remark}{\rm
As we have already seen, the technique of \'etale coverings requires
 special care in the case of positive characteristics.
For example, let $\EE = \EE_1$ be a Weierstra\ss{} nodal curve and
$\pi_2: \EE_2  \lar  \EE_1$ an \'etale covering of degree $2$. Then the vector bundle
$\pi_{2*}(\kO)$ corresponds to  the triple
$(\widetilde\kO^2, \kk^2(s), \tilde{i})$, where
$\tilde{i}$ is given by matrices
$$
i(0) =
\setlength\arrayrulewidth{0.1mm}
\setlength\doublerulesep{0.1mm}\doublerulesepcolor{black}
\begin{array}{
@{}c@{}     c       @{}c@{}     c       @{}c@{}  }
\hline\hline
\fvr &\mk{1} &   &\mk{0} &\fvr \\
\fvr &\mk{0} &  &\mk{1} &\fvr \\
\hline\hline
\end{array}
~~\hbox{and}~~
i(\infty) =
\begin{array}{
@{}c@{}     c       @{}c@{}     c       @{}c@{}  }
\hline\hline
\fvr &\mk{0} &   &\mk{1} &\fvr \\
\fvr &\mk{1} &  &\mk{0} &\fvr \\
\hline\hline
\end{array}$$
Then for $char(\kk) \ne 2$ we have
$$
\setlength\arrayrulewidth{0.1mm}
\setlength\doublerulesep{0.1mm}\doublerulesepcolor{black}
\begin{array}{
@{}c@{}     c       @{}c@{}     c       @{}c@{}  }
\hline\hline
\fvr &\mk{0} &   &\mk{1} &\fvr \\
\fvr &\mk{1} &  &\mk{0} &\fvr \\
\hline\hline
\end{array}
\sim
\setlength\arrayrulewidth{0.1mm}
\setlength\doublerulesep{0.1mm}\doublerulesepcolor{black}
\begin{array}{
@{}c@{}     c       @{}c@{}     c       @{}c@{}  }
\hline\hline
\fvr &\mk{1} &   &\mk{0} &\fvr \\
\fvr &\mk{0} &  &\mk{-1} &\fvr \\
\hline\hline
\end{array}$$
and $\pi_{2*}\kO \cong \kO \oplus \kB(0,1, -1)$. However, for $char(\kk) =  2$
$$
\setlength\arrayrulewidth{0.1mm}
\setlength\doublerulesep{0.1mm}\doublerulesepcolor{black}
\begin{array}{
@{}c@{}     c       @{}c@{}     c       @{}c@{}  }
\hline\hline
\fvr &\mk{0} &   &\mk{1} &\fvr \\
\fvr &\mk{1} &  &\mk{0} &\fvr \\
\hline\hline
\end{array}
\sim
\begin{array}{
@{}c@{}     c       @{}c@{}     c       @{}c@{}  }
\hline\hline
\fvr &\mk{1} &   &\mk{1} &\fvr \\
\fvr &\mk{0} &  &\mk{1} &\fvr \\
\hline\hline
\end{array}
$$
and $\pi_{2*}\kO \cong \kF_2$.}
\end{remark}

From Theorem \ref{thm:etale} one can derive  formulas for the cohomology groups of
indecomposable torsion free sheaves, a formula for  the dual
of an indecomposable torsion free sheaf and rules for the computation of
the direct sum decomposition of two indecomposable vector bundles. This is what we are going to
describe now.

\begin{lemma}[\cite{BDG, BurbanKreussler1}]
  If $\boldsymbol{d}=(d_{1},\ldots,d_{rn})\in \mathbb{Z}^{rn},
\boldsymbol{e}= (e_1, e_2,\dots, e_k),$
$  \lambda\in\boldsymbol{k}^{*}, m\ge 1, 1\le f \le n$, we have:
  \begin{itemize}
  \item[(i)] $\mathcal{B}(\boldsymbol{d},m,\lambda)^{\vee} \cong
              \mathcal{B}(-\boldsymbol{d},m,\lambda^{-1})$
  \item[(ii)] $\mathcal{S}(\boldsymbol{e}, f)^{\vee} \cong
               \mathcal{S}(\boldsymbol{\kappa - e}, f)$ with
              $\boldsymbol{\kappa}=
              \begin{cases}
                (-1,0,\ldots,0,-1) &\text{if } k \ge 2\\
                -2             &\text{if } k = 1.
              \end{cases}$
  \end{itemize}
\end{lemma}

\noindent
\emph{Proof}.
  If $f: \XX \rightarrow \EE$ is a finite morphism,
$\kF$ a  coherent sheaf
  on $\XX$ and $\kG$ a locally free sheaf on $\EE$ then there is a natural isomorphism of
  $f_{\ast}\mathcal{O}_{\XX}$--modules
  $$f_{\ast}\mathcal{H}om_{\XX}(\kF,f^{!}\kG) \cong
            \mathcal{H}om_{\EE}(f_{\ast}\kF,\kG).$$
    Recall that $f^{!}\omega_{\EE}$ is a dualizing sheaf on $\XX$ if
    $\omega_{\EE}$ is one on $\EE$. In our situation
    $\omega_{\EE} \cong \mathcal{O}_{\EE}$ and we obtain
    an isomorphism
    $$f_{\ast}\mathcal{H}om_{\XX}(\kF,\omega_{\XX}) \cong
    \mathcal{H}om_{\EE}(f_{\ast}\kF,\mathcal{O}_{\EE})
    \cong (f_{\ast}\kF)^{\vee}.$$
    To show (i), we consider $\XX=\EE_{n}$ and $f=\pi_{n}$. The claim
    follows now from $\omega_{\EE_{n}} \cong
    \mathcal{O}_{\EE_{n}}$, $\mathcal{F}_{m}^{\vee}
    \cong \mathcal{F}_{m}$ and
    $\mathcal{L}(\boldsymbol{d},\lambda)^{\vee} \cong
    \mathcal{L}(-\boldsymbol{d},\lambda^{-1})$ on $\EE_{n}$.

    For the proof of (ii) we let $\XX=\mathbb{I}_{k}$ and $f=p_{k}$. Now
    $\omega_{\mathbb{I}_{k}} \cong \mathcal{L}(\boldsymbol{\kappa})$ and
    the result follows from $\mathcal{L}(\boldsymbol{d})^{\vee} \cong
    \mathcal{L}(-\boldsymbol{d})$ on $\mathbb{I}_{k}$.

\medskip
\medskip
Using the description of indecomposable vector bundles  via \'etale coverings
it is not difficult to compute their cohomology.

\begin{lemma}[\cite{DGK}]
There is the following formula for the cohomology of
indecomposable vector bundles
$$
dim_{\kk} H^{0}(\kB(\boldsymbol{d},m,\lambda)) = m\bigl(\sum\limits_{i=1}^{rn} (d_i + 1)^{+} -
\theta({\boldsymbol d})\bigr) + \delta({\boldsymbol  d},\lambda)
$$
and
$$
dim_{\kk} H^{1}(\kB(\boldsymbol{d},m,\lambda)) = rm -
dim_{\kk} H^{0}(\kB(\boldsymbol{d},m,\lambda),
$$
where $\delta({\boldsymbol  d},\lambda)=1$ if ${\mathbf d}=(0,\dots,0)$, $\lambda=1$ and
$0$ otherwise; $k^{+} = k$ if $k > 0$ and zero otherwise.
The number $\theta({\boldsymbol d})$ is defined as follows:
call a subsequence ${\mathbf p}=(d_{k+1},\dots, d_{k+l})$, where $0\le k < rn$ and
$1\le l \le rn$  a \emph{positive part} of ${\boldsymbol  d}$ if all $d_{k+j}\ge 0$ and
either $l=rn$ or both $d_{k}<0$ and $d_{k+l+1}<0$. For such a positive part put
$\theta({\boldsymbol  p})=l$ if either $l=rs$ or ${\mathbf p}=(0,\dots,0)$ and
$\theta({\boldsymbol p})=l+1$
otherwise. Then $\theta({\boldsymbol  d})=\sum \theta({\boldsymbol  p})$, where we take a sum over
all positive subparts of ${\boldsymbol d}$.
\end{lemma}

In order to compute the tensor product of two indecomposable vector bundles we shall need the following lemma.

\begin{lemma}
Let
$\EE$
be a cycle of projective lines,
$\pi_i: \EE_i \lar \EE$ two \'etale coverings  $i = 1,2$ and
 $\kE_i$ a vector bundle on $\EE_i$. Let
$\EE'$ be the  fiber product of $\EE_1$ and $\EE_2$ over $\EE$:
$$
\begin{CD} \EE' @>p_{1}>> \EE_1 \\ @VVp_{2}V
          @VV\pi_{1}V \\ \EE_2 @>\pi_{2}>> \EE.
\end{CD}
$$
Denote by
$\widetilde{\pi}: \EE' \lar \EE$  the composition  $\pi_1 p_1$,  then
$$
  \pi_{1*}(\kE_{1})\otimes \pi_{2*}(\kE_{2}) \cong
\widetilde{\pi}_{*}\bigl(p_{1}^{*}(\kE_{1})\otimes p_{2}^{*}(\kE_{2})\bigr).
$$
\end{lemma}

\noindent
\emph{Proof}.
By the base change and projection formula
\begin{eqnarray*}
\widetilde{\pi}_{*}\bigl(p_{1}^{*}(\kE_{1})\otimes p_{2}^{*}(\kE_{2})\bigr) \cong
\pi_{2*} p_{2*} \bigl(p_{1}^{*}(\kE_{1}) \otimes p_{2}^{*}(\kE_{2})\bigr) &\\
\cong
\pi_{2*}\bigl(p_{2*} p_{1}^{*} (\kE_{1}) \otimes \kE_{2}\bigr) \cong
 \pi_{2*}\bigl(\pi_{2}^{*} \pi_{1*} (\kE_{1}) \otimes \kE_{2}\bigr) \cong&
\pi_{1*}(\kE_{1}) \otimes \pi_{2*}(\kE_{2}).
\end{eqnarray*}

\medskip

The following proposition describes the fiber product of two \'etale coverings of a
given cycle of
projective lines.

\begin{proposition}[\cite{Burban}]
Let  $\EE_n$ be a cycle of $n$ projective lines and
 $\pi_{i}: \EE_{d_i n} \lar \EE_n$
two \'etale covering of degree $d_i$, $i=1,2$.  Choose a labeling
of  the irreducible components  of  each cycle
$\EE_{d_1 n}$ and $\EE_{d_2 n}$ by consecutive non-negative integers
$0,1,2, \dots, $ such
that  $\pi_i$  maps the zero component into a zero component, $i=1,2$.
Let $d=g.c.d.(d_{1},d_{2})$ be the greatest  common divisor and
 $D=[d_{1},d_{2}]$ the smallest common multiple of  $d_1$ and $d_2$. Then the fiber product is
$\widetilde{\EE} = \coprod\limits_{i=1}^d \EE_{Dn}^{(i)}$ and
$p_{1}: \EE_{Dn}^{(i)}\lar \EE_{d_1 n}$  is the \'etale  covering determined by the assumption
that it maps the
$i$-th component of   $\EE_{Dn}^{(i)}$ to the 0-th component of
$\EE_{d_{1} n}$.  The second morphism
$p_{2}:  \EE_{Dn}^{(i)}\lar \EE_{d_{2} n}$ is  the \'etale covering,  mapping
the  zero component to the zero component.
\end{proposition}

These properties  allow
to describe a decomposition of the tensor product of any two indecomposable
vector bundles into a direct sum of indecomposable ones. In particular, in  the case of a nodal
Weierstra\ss{} curve we get  the following concrete algorithm, obtained for the
first time in  \cite{Yudin}.

\begin{theorem}[\cite{Yudin, Burban}]
Let $\EE$ be a Weierstra\ss{} nodal curve over an algebraically closed field $\kk$
of characteristics zero,
$\kB({\boldsymbol  d},1,\lambda)$ and $\kB({\boldsymbol  e},1,\mu)$  two
vector bundles on $\EE$ of rank $k$ and $l$ respectively,
$\boldsymbol{d} = d_{1}d_{2}\dots d_{k}$ and $\boldsymbol{e} =
e_{1}e_{2}\dots e_{l}$. Let $D$ be the smallest common multiple
and $d$ the greatest  common divisor of
$k$ and  $l$. Consider  $d$ sequences
$$
\begin{array}{l}
  \boldsymbol{f}_{1} = d_{1}+e_{1},d_{2}+e_{2},\dots,d_{k}+e_{l},\\
  \boldsymbol{f}_{2} = d_{1}+e_{2},d_{2}+e_{3},\dots,d_{k}+e_{1},\\
\vdots \\
\boldsymbol{f}_{d} = d_{1}+e_{d},d_{2}+e_{d+1},\dots,d_{k}+e_{d-1},\\
\end{array}
$$
of length $D$. Then the following decomposition holds:
$$
  \kB(\boldsymbol{d},1,\lambda)\otimes \kB(\boldsymbol{e},1,\mu) \cong \bigoplus_{i=1}^{d}
\kB(\boldsymbol{f_{i}},1,\lambda^{\frac{l}{d}} \mu^{\frac{k}{d}}).
$$
If some $\boldsymbol{f_{i}}$  is periodic, then we use  the  isomorphism
$$
  \kB(\boldsymbol{g}^{l},1,\lambda) = \bigoplus_{i=1}^{l}
\kB(\boldsymbol{g},1,\xi^{i}\sqrt[l]{\lambda}),$$
where  $\boldsymbol{g}^{l} = \underbrace{\boldsymbol{g}\boldsymbol{g}\dots
\boldsymbol{g}}_{l}$
and   $\xi$ a primitive $l$--th root   of $1$.
\end{theorem}

\medskip

Even possessing  a complete classification of indecomposable torsion free sheaves on a Weierstra\ss{}
nodal curve,  an exact  description of stable vector bundles is a non-trivial problem.
It  can be shown by many methods that for a pair of coprime integers
$(r,d) \in \mathbb{Z}^2$, $r>0$
the moduli space of stable vector bundles of rank $r$ and degree
$d$ is $\kk^*$, see for example \cite{BurbanKreussler3}. However, for applications
it is  important to have a description of stable vector bundles via \'etale coverings.
In order to get such a classification note the following useful fact.

\begin{lemma}[\cite{Burban}]
\label{lem:stabsimpl}
Let $\EE$ be an irreducible Weierstra\ss{} curve. Then a coherent sheaf $\kF$ on $\EE$
is stable if and only if it is simple i.e. $\End(\kF)=\kk.$
\end{lemma}

In general, for irreducible curves
 stability implies simplicity, but  in the case of irreducible curves of arithmetic genus one
 both conditions are equivalent.
Then one can prove the following theorem:
%Using the whole machinery we developed for the study of vector bundles on cycles of projective lines one can show the following theorem.

\begin{theorem}[\cite{Burban}]
\label{thmBurbStable}
Let $\EE$ be a nodal Weierstra\ss{} curve and  $\kE$  a stable vector bundle on $\EE$ of
rank $r$ and degree $d$, $0 < d < r$.  Then $g.c.d.(r,d) = 1$, $\kE \cong
\kB(\boldsymbol{d},1,\lambda)$ and $\boldsymbol{d}$ can be obtained by
the following  algorithm.

\begin{enumerate}
\item Let $y = min(d, r-d)$, $x = max(d, r-d)$. If $x = y$, then
$\boldsymbol{d} = (0,1)$.
Assume now $x >y$. Consider the triple $(x, y, x+y)$  and write
$x+y = (k+1)y + s$, where $0 < s < y$ and  $k \ge 1$.  If $s > y-s$ then replace
$(x, y, x+y)$ by $(s, y-s, y)$ and  say
say that $(x,y,x+y)$ is obtained from $(s, y-s, y)$ by the blow-up of type
$(A,k)$. If $s <  y-s$ then replace  $(x, y, x+y)$ by $(y-s, s, y)$ and say that
$(x,y,x+y)$ is obtained from $(y-s, s, y)$ by the
blow-up of type   $(B,k)$.
\item Repeat this algorithm until we get the triple $(p, 1, p+1)$.
Consider the sequence of reductions $(x, y, x+y) =
(x_0, y_0, x_0 + y_0)  \stackrel{(C_1, k_1)}\lar (x_1, y_1, x_1 + y_1) \lar
\dots  \stackrel{(C_n, k_n)}\lar (x_n, y_n, x_n + y_n) = (p, 1, p+1)$, where
$C_i \in \{A, B\}$ and $k_i \ge  1$ for $1 \le i \le n$.
\end{enumerate}

Now we can recover the vector $\boldsymbol{d}:$
\begin{enumerate}
\item Start with  sequence $\underbrace{\alpha,\alpha,\dots,\alpha}_{p\, times},\beta$,
which corresponds
to the triple $(x_n, y_n, \\ x_n + y_n) = (p, 1, p+1)$. If $C_n = A$ then
replace each letter
$\alpha$  by the block
$\underbrace{\alpha,\alpha,\dots,\alpha}_{\mbox{\tiny $k_n+1$}}$
and  each  $\beta$ by the block
$\underbrace{\alpha,\alpha,\dots, \alpha}_{\mbox{\tiny $k_n$}}$. Between these new blocks insert
the letter  $\beta$. If  $C_n = B$ then
replace each letter    $\alpha$  by the block
$\underbrace{\alpha,\alpha,\dots,\alpha}_{\mbox{\tiny $k_n$}}$
and each letter $\beta$ by the block
$\underbrace{\alpha,\alpha,\dots,\alpha}_{\mbox{\tiny $k_n+1$}}$. Between these  new blocks
insert   the letter  $\beta$ again. We have got a new sequence of letters $\alpha$ and $\beta$ of
total length $x_{n-1} + y_{n-1}$ with $x_{n-1}$ letters $\alpha$ and $y_{n-1}$ letters $\beta$.
\item Proceed inductively until we get  a sequence of length $r$ with $max(d, r-d)$
letters  $\alpha$ and $min(d, r-d)$ letters $\beta$.
\item If $d > r-d$ then
replace each letter   $\alpha$ by $1$ and each letter $\beta$ by $0$. In case  $d \le r-d$
replace   $\alpha$ by $0$ and $\beta$ by $1$.  The
resulted sequence is the vector $\boldsymbol{d}$ we are
looking for.
\end{enumerate}
\end{theorem}

\begin{example}
{\rm Let rank $r = 19$ and degree $11$. The sequence of reductions is
$$
(11,8,19) \stackrel{(B,1)}\lar (5,3,8) \stackrel{(A,1)}\lar
(2,1,3).
$$
Using the algorithm we get the sequence of blowing-ups
\begin{align*}
\alpha,\alpha,\beta  &\stackrel{(A,1)} \lar  \alpha,\alpha,\beta,\alpha,\alpha,\beta,\alpha,\beta &\\
\alpha,\alpha,\beta,\alpha,\alpha,\beta,\alpha,\beta
&
\stackrel{(B,1)}\lar \alpha,\beta,\alpha,\beta,\alpha,\alpha,\beta,\alpha,\beta,
\alpha,\beta,\alpha,\alpha,\beta,\alpha,\beta,\alpha,\alpha,\beta
\end{align*}
and hence
$$
\boldsymbol{d} =  (1,0,1,0,1,1,0,1,0,1,0,1,1,0,1,0,1,1,0).
$$
}
\end{example}

This result was generalized by Mozgovoy \cite{Mozgovoy} to get a recursive description of
semi-stable torsion free sheaves of arbitrary slope.

\medskip
\medskip

Note the following important difference between  smooth and singular curves of arithmetic
genus one. In the smooth case any indecomposable coherent sheaf is either
locally free or torsion free and is  automatically
semi-stable. This is no longer true for singular curves, in particular, in that case
there are indecomposable
coherent sheaves which are neither torsion nor torsion free.

\begin{example}
{\rm
Let $\EE$ be a nodal Weierstra\ss{} curve, $s$ its singular point, $n: \PP^1 \lar \EE$
the normalization map. Then $$\Ext^1(n_*(\kO_{\PP^1}), \kk(s)) =
H^0({\mathcal Ext}^1(n_*(\kO_{\PP^1}), \kk(s))) = \kk^2.$$
Let $w \in \Ext^1(n_*(\kO_{\PP^1}), \kk(s))$
be a non-zero element and
$$
0 \lar \kk(s) \stackrel{i}\lar \kF \stackrel{p}\lar n_*(\kO_{\PP^1}) \lar 0
$$
the corresponding extension. Then $\kF$ is an indecomposable coherent sheaf which is neither
torsion nor torsion free.  To see that $\kF$ is indecomposable assume
$\kF \cong \kF' \oplus \kF''$. Then one of  its   direct summands, say $\kF'$ is a torsion sheaf.
Since $\Hom(\kF', n_*(\kO_{\PP^1})) = 0$, $\kF'$ belongs to the kernel $ker(p)$
and hence is isomorphic to $\kk(s)$. Therefore  the map $i$ has   a left inverse,
hence $w = 0,$ and that is a contradiction.
}
\end{example}

\begin{proposition}[\cite{BurbanKreussler3}]\label{prop:extreme}
  Let $\EE$ be a singular Weierstra\ss{} curve and
 $\kF \in \Coh(\EE)$  an indecomposable coherent sheaf which is not
  semi-stable. Then, all Harder-Narasimhan factors of $\kF$ are direct sums
of semi-stable sheaves of infinite homological dimension.
\end{proposition}

\noindent
\emph{Proof}.
Let
 $0\subset \mathcal{F}_{n} \subset \ldots \subset \mathcal{F}_{1} \subset
\mathcal{F}_{0} = \mathcal{F}$ be the Harder-Narasimhan filtration of $\kF$
with semi-stable factors $\mathcal{A}_{\nu} := \mathcal{F}_{\nu}/\mathcal{F}_{\nu+1}$ of
decreasing slopes $\mu(\mathcal{A}_{n})>
\mu(\mathcal{A}_{n-1}) > \ldots > \mu(\mathcal{A}_{0})$.

Assume $\kA_\nu \cong \kA_{\nu}' \oplus \kA_{\nu}''$ and $\kA_{\nu}'$ has finite
global dimension.
Since $\kF_{\nu+1}$ is filtered by semi-stable sheaves $\kF_\mu$ for  $\mu > \nu$ and
$\Hom(\kA_\mu, \kA_{\nu}') = 0$, we get
$\Ext^1(\kA_{\nu}', \kF_{\nu+1}) \cong \Hom(\kF_{\nu+1}, \kA_{\nu}')^* = 0$.
Therefore $\kF_\nu$ contains $\kA_{\nu}'$ as a direct summand: $\kF_\nu \cong
\kF_{\nu}' \oplus \kA_\nu$. From the exact sequence
$$
0 \lar \kF_{\nu}' \oplus \kA_{\nu}'  \lar \kF_{\nu - 1} \lar \kA_{\nu-1} \lar 0
$$
and the isomorphism $\Ext^1(\kA_{\nu-1}, \kA_{\nu}') \cong \Hom(\kA_{\nu}', \kA_{\nu-1})^* = 0$
we conclude that $\kF_{\nu - 1}$ contains $\kA_{\nu}'$ as a direct summand as well.
Proceeding inductively we obtain that $\kF$ itself contains $\kA_\nu'$ as a direct summand, a
contradiction.

We see that a difference between the combinatorics of indecomposable coherent sheaves
on smooth and singular Weierstra\ss{} curves is due to the existence of
semi-stable sheaves of infinite global dimension together with the failure of the
Serre duality on singular curves.  In order to classify indecomposable coherent sheaves it is convenient
to
consider a more general problem: the  description of indecomposable objects of the derived
category
$D^-(\Coh(\EE))$.  It turns out that the last problem is again tame and can be solved using
the technique of representations of bunches of chains, see \cite{BurbanDrozd1} for the details.

\def\vb{\VB} %\mathrm{VB}}
\def \oo {\kO}
\def \ko {{\bf k}}
\def \en  {\widetilde{\EE}}

\def \on  {\widetilde{\oo}}
\def \aa  {{\kA}}
\def \an  {\widetilde{\kA}}
\def \jj  {{\kJ}}
\def \ff  {\kF}
\def \fn  {\widetilde{\kF}}
\def \bn  {\widetilde{\kE}}
\def \mo    {\rightarrow}
\def \hh  {{H}}

%to matrix
%\def \vr    {\vrule}
%\def \br    {\bold{|}}
%\newcommand{\br}{\vrule\makebox[0.1pt][c]{}\vrule}

\subsection{Vector bundles on a cuspidal cubic curve}

As we have mentioned in the introduction, the category of vector bundles
on a curve of arithmetic genus one, different from a cycle of projective lines, is
vector bundle wild, see also Corollary \ref{cor:cuspiswild}.
Nevertheless,
 if we restrict ourselves to the subcategory of simple
 vector bundles
$\vb_s,$ or even to the subcategory of simple torsion free sheaves $\TF_s$,
then the classification problem becomes tame again and, moreover, the combinatorics of
the answer resembles %very much
the case of smooth and nodal Weierstra\ss{} curves (see
Theorem \ref{propCoh}  and  Theorem \ref{thmBurbStable}).

 \begin{theorem}\label{01}
  Let $\EE$ be a cuspidal cubic curve over an algebraically closed field $\ko.$ Then
 \begin{enumerate}
\item
 the rank $r$ and the degree $d$ of a simple torsion free sheaf $\kF$ over $\EE$ are coprime;
 \item
  for every pair $(r,d)$ of coprime integers with positive $r,$
the isomorphism classes of
simple vector bundles $\kE\in\VB_s(r,d)$ are parametrized by $\mA^1$
and
there is a
unique simple torsion free but not locally free sheaf $\ff$ of rank $r$ and degree $d.$
%up to an isomorphism.
\end{enumerate}

%\begin{remark}
\rm Note that $\mA^1\cong\EE_{reg}$
 is isomorphic to the Picard group $\Pic^{\circ}(\EE).$
\rm
It can be shown that for a pair of coprime integers $r>0$ and $d$
the moduli  space of $\TF_s(r,d)$ is isomorph to $\EE,$ moreover,
vector bundles $\kE$ correspond to nonsingular points of $\EE$   and
the unique torsion free but not locally free sheaf $\ff$ corresponds to
the singular point $s$.
%\end{remark}
%\cite[Example II.6.11.4]{ha}.
\end{theorem}

\noindent
\emph{Sketch of proof}.
Let $\EE$ be a cuspidal cubic curve given by the equation  $x^3-y^2z=0.$
Choose coordinates $(z_0: z_1)$ on the normalization $\widetilde\EE \cong \mP^1
\stackrel{p}\lar \EE$ such that the preimage of
the singular point $s = (0:0:1)$ of $\EE$  is $(0:1)$. Let $U = \{(z_0:z_1)| z_1 \ne 0\}$
be an affine neighborhood of $(0:1)$ and $z = z_0/z_1$.
In the notations of Section \ref{subsecTr} we have: $\aa\cong\ko(s)$ and
$\an\cong \big(\ko[\varepsilon]/\varepsilon^2\big)(s).$

Let $\kF$ be a torsion free sheaf of rank $r$ on $\EE$ and
$(\widetilde\kF, \kM, \tilde{i})$ be the  corresponding triple. Then, as in the case of a nodal
rational curve, we have
\begin{itemize}
\item
 a splitting $\widetilde\kF \cong
\bigoplus\limits_{n \in \mathbb{Z}} \kO(n)^{r_n}$, with $\sum_{n\in\mZ} r_n=r;$
\item an isomorphism $\kM \cong \kA^t,$ for some $t\geq r,$ and
 $t=r$ if and only if $\ff$ is a vector bundle;
\item  an epimorphism of
$\widetilde\kA$--modules $\tilde{i}: \widetilde\kM \otimes_{\kA} \widetilde\kA \lar
\widetilde\kF\otimes_{\widetilde\kO} \widetilde\kA,$ which is an isomorphism
if and only if $\ff$ is a vector bundle.
\end{itemize}
In order to write $\tilde{i}$ in matrix form remember that we
identify $\widetilde{\kF}\otimes_{\widetilde\kO} \widetilde\kA$ with
$$p_*\big(p^*(\kF) \otimes_{\kO_{\mP^1}} \kO_{\mP^1}/\kI\big),  \hbox{~where~}~
\kI = \kI^2_{(0:1)}$$
is the ideal sheaf  of the scheme-theoretic preimage of $s$.
We choose a basis of
$\kM \cong \ko(s)^r$ and fix the trivializations
$$
\kO_{\mP^1}(n) \otimes  \kO_{\mP^1}/\kI \lar
\big(\ko[\varepsilon]/\varepsilon^2\big)(s)
$$
given by the map  $\zeta\otimes 1 \mapsto pr(\frac{\zeta}{z_1^n})$ for a local section $\zeta$
of $\kO(n)$ on an open set $V$ containing $(0:1),$
where $pr: \ko[V] \lar \ko[\varepsilon]/\varepsilon^2$
be the map induced by  $\ko[z] \lar \ko[\varepsilon]/\varepsilon^2$,
$z \mapsto \varepsilon$.
Using these choices we may write $\tilde{i} = i(0) + \varepsilon i_\varepsilon(0),$ where
both $i(0)$ and $i_\varepsilon(0)$ are square $r\times r$ matrices. Since by Theorem
\ref{thm:triples}
the isomorphism classes of triples are in bijection with the isomorphism classes of
torsion free sheaves,
we have
to study the action of automorphisms  of $(\widetilde\kF, \kM, \tilde{i})$ on the
matrices $i(0)$ and $i_\varepsilon(0)$.
The condition for $\tilde{i}$ to be surjective is
equivalent to the surjectivity
%surj
 of $i(0)$.
Similarly, for vector bundles we have  that
$\tilde{i}$ is invertible if and only if
 $i(0)$ is invertible.

If we have a morphism  $\kO(n) \lar \kO(m)$ given by a homogeneous form
$Q(z_0,z_1)$ of degree $m-n$ then the induced map
$\kO(n) \otimes \kO/\kI  \lar \kO(m) \otimes \kO/\kI$ is given by the map
$pr(Q(z_0, z_1)/z_1^{m-n}) = Q(0:1) + \varepsilon \frac{d Q}{d z_0}(0:1)$.

Moreover,
for any endomorphism $(F,f)$ of the triple
$(\widetilde\kF, \kM, \tilde{i})$ the induced map
 $\bar{F}:\fn \otimes \kO/\kI  \lar \fn \otimes \kO/\kI$
has the form $\bar{F}  = F(0:1)  + \varepsilon \frac{d F}{d z_0}(0:1).$
If $(F,f)$ is an automorphism  then $\bar{F}\in GL_r(\ko[\varepsilon]/\varepsilon^2)$ and
 the transformation rule
$\bar{F}\tilde{i} = \tilde{i'}{f}$  in matrix form reads
$$
\begin{array}{ccl}
i'(0) &=&  F(0:1) i(0) f^{-1} \\
i'_\varepsilon (0)& = & \frac{d F}{d z_0}(0:1)i(0) f^{-1} + F(0:1) i_\varepsilon(0) f^{-1}.
\end{array}
$$
As a result, the  matrix problem is as follows:
we have two matrices $i(0)$ and $i_\varepsilon(0)$
with $r$ rows and $t$ columns, and
 $rank(i(0))=r.$
In the case of a vector bundle $i(0)$ and $i_\varepsilon(0)$ are square matrices
and $i(0)$ is invertible.
The matrices $i(0)$ and $i_\varepsilon(0)$
 are divided into horizontal blocks labelled by integers called weights. Any two blocks
of $i(0)$ and $i_\varepsilon(0)$ marked by the same label
are called {\it conjugated} and have the same number of rows.

\begin{figure}[ht]
\hspace{2cm}
\begin{minipage}[b]{4cm}
\begin{align*}
\setlength\arrayrulewidth{0.3mm}
\begin{array}{ @{} c@{} c      @{}c@{} l }
\cline{1-2}
\fvr & %\hspace{3cm} d%\phantom{abrakadabradabra}
{~\mk{}\mk{}  {\vdots}  \mk{}\mk{}~}
& \fvr\\
%\cline{1-2}
%\fvr &  {\vdots} & \fvr\\
\cline{1-2}
\fvr &{\scriptstyle n-1} & \fvr\\
\cline{1-2}
\fvr & {\scriptstyle n} & \fvr & \} {\scriptstyle r_n} \\
\cline{1-2}
\fvr &_{\scriptstyle n+1} & \fvr\\
\cline{1-2}
\fvr & ^{\vdots} & \fvr\\
\cline{1-2}
%\fvr & & \fvr\\
%\cline{1-2}
\end{array}
\end{align*}
\center {$i(0)$}
%\caption*{c)}
\end{minipage}
\begin{minipage}[b]{4cm}
\begin{align*}
\setlength\arrayrulewidth{0.3mm}
\begin{array}{ @{} c@{} c      @{}c@{} l }
\cline{1-2}
\fvr & %\hspace{3cm} d%\phantom{abrakadabradabra}
{~\mk{}\mk{}  {\vdots}  \mk{}\mk{}~}
& \fvr\\
%\cline{1-2}
%\fvr &  {\vdots} & \fvr\\
\cline{1-2}
\fvr &{\scriptstyle n-1} & \fvr\\
\cline{1-2}
\fvr & {\scriptstyle n} & \fvr & \} {\scriptstyle r_n} \\
\cline{1-2}
\fvr &_{\scriptstyle n+1} & \fvr\\
\cline{1-2}
\fvr & ^{\vdots} & \fvr\\
\cline{1-2}
%\fvr & & \fvr\\
%\cline{1-2}
\end{array}
\end{align*}
\center {$i_{\varepsilon}(0)\phantom{XX}$}
%\caption*{c)}
\end{minipage}
\end{figure}

The permitted transformations are listed below:
\begin{enumerate}
\item We can simultaneously do any elementary transformations of columns of
$i(0)$
and of $i_\varepsilon(0)$.
\item We can simultaneously perform any invertible elementary transformations
of rows of $i(0)$ and $i_\varepsilon(0)$
inside of any two  conjugated horizontal blocks.
\item We  can  add
a scalar multiple of any row with lower weight to any row with higher weight
simultaneously in $i(0)$ and $i_\varepsilon(0).$
\item We can
add a row  of $i(0)$ with a lower weight
to any row of
$i_\varepsilon(0)$ with a
higher weight.
\end{enumerate}
%(Compere this with the case of nodal Weierstra\ss{} curve!)
%???
This matrix problem turns to be wild, see corollary \ref{cor:cuspiswild}.
However, the simplicity  condition of a triple
 $(\widetilde{\kF}, {\kM}, \tilde{i})$ implies additional restrictions, which make
 the problem tame. First note that
if  $\fn$ contains $\on(c)\oplus \on(d)$ with $d>c+1$ as a direct summand, then
the pair $(F,f):=(z_0^{d-c}, 0)$ defines a non-scalar endomorphism
of the triple $(\fn, \kM, \tilde{i}),$ which, therefore, can not be simple.
Thus, for a  simple torsion free sheaf $\kF$ we may  assume
$\fn\cong\on(c)^{r_1}\oplus \on(c+1)^{r_2}$ for some $c\in\mZ$ and
the matrix $\tilde{i}$
consists of two horizontal blocks.

We consider the case of vector bundles first.
Although the case of
torsion free but not locally free sheaves
is similar, it should
be considered separately in order to make the presentation
clearer.

%Let us concentrate on the case of vector bundles first.
As was mentioned above, if $\ff$ is a vector bundle  then
$\tilde{i}$ is an isomorphism and by transformations 1 and 2
the matrix $i(0)$ can be reduced to the identity matrix.
Moreover,
 by applying transformation 4 we can make the left lower block of $i_\varepsilon(0)$ zero,
 as indicated below:
%\begin{equation}
$$
\setlength\arrayrulewidth{0.1mm}
\setlength\doublerulesep{0.1mm}\doublerulesepcolor{black}
i(0)=
\begin{array}{
@{}c@{}     c       @{}c@{}     c       @{}c@{}}
\hline\hline
%    & \hr && \hr \\
\fvr  & {\mk{I_{r_1}}}  &\fvr&    {\mk{}}   & \fvr     \\
\hline\hline
\fvr&       \mk{}        &\fvr &        {\mk{I_{r_2}}} & \fvr \\
%\cline{3-4}
\hline\hline
\end{array}
~~~\hbox{and}~~~~
i_\varepsilon(0)=
\begin{array}{
@{}c@{}     c       @{}c@{}     c       @{}c@{}}
\hline\hline
%    & \hr && \hr \\
\fvr  & {\mk{B_{1}}}  &\fvr&    {\mk{B_{12}}}   & \fvr     \\
\hline\hline
\fvr&       \mk{}        &\fvr &        {\mk{B_{2}}} & \fvr \\
%\cline{3-4}
\hline\hline
\end{array}.
\label{eqi(0)}
\eqno(*)
$$
%\end{equation}
Here $I_n$ denotes  the identity matrix of size $n,$
an empty space stands for a zero block
and $B_1, B_{12}, B_2$ denote
nonreduced blocks.

Thus we can assume that $i(0)$ is the identity matrix and concentrate on the matrix
$i_\varepsilon(0)$, taking into account only those transformations, which leave
$i(0)$ unchanged. Then we obtain the category of block matrices
$$\BM=\bigcup\limits_{(r_1,r_2)}\BM(r_1,r_2).$$
Objects of $\BM(r_1,r_2)$
are matrices of the form $i_{\varepsilon}(0)$ in formula (*), i.e.
upper triangular
block matrices $B$  consisting of the blocks $(B_1,B_{12},B_2),$
where $(B_1,B_2)$ are square matrices of sizes $r_1$ and $r_2$ respectively.
Morphisms  $C: B\mo B^{\prime}$  are given by lower triangular block matrices:
%consisting of blocks $(C_1, C_{21}, C_2 )$
$$
\setlength\arrayrulewidth{0.1mm}
\setlength\doublerulesep{0.1mm}\doublerulesepcolor{black}
C=
\begin{array}{
@{}c@{}     c       @{}c@{}     c       @{}c@{}}
\hline\hline
%    & \hr && \hr \\
\fvr  & {\mk{C_{1}}}  &\fvr&    \mk{}   & \fvr     \\
\hline\hline
\fvr&       {\mk{C_{21}}}        &\fvr &        {\mk{C_{2}}} & \fvr \\
%\cline{3-4}
\hline\hline
\end{array}
%\begin{array}{@{}c@{}     c}
%\left\} r_1\makebox[1pt][c]{\rule{0pt}{10pt}}\right.\\
%\left\} r_2\mk{}\right.
%\end{array}
$$
with block sizes $(r_1, r_2)$ and
satisfying equations $CB=B^{\prime}C.$ In term of blocks this equation can be written as:
$$
 \label{eqBlock}
  \begin{array}{rcl}
   C_1B_1&=&B^{\prime}_1C_1+B^{\prime}_{12} C_{21},\\
  C_1B_{12}&=&B^{\prime}_{12}C_2,\\
  C_2B_2+C_{21}B_{12}&=&B^{\prime}_2C_2.
 \end{array}
\eqno(**)
$$

Two matrices
$B$ and $B^{\prime}$
are called {\it equivalent} (i.e. correspond to isomorphic vector bundles) if there
is a non-degenerate  morphism $C: B\mo B^{\prime},$ i.e. if $B^{\prime}=CBC^{-1}.$
In terms of transformations this means: we can  add a row $k$ with lower weight
to a row $j$ with higher weight and simultaneously add the column $j$ to the column $k.$
A matrix $B\in \BM(r_1, r_2)$ is called {\it simple } if any endomorphism $C: B\mo B$
is scalar.
Obviously, simplicity is a property defined on equivalence classes.
The full subcategory $\BM(r_1, r_2)$ consisting of simple objects $B$
is denoted by $\BM_s(r_1, r_2).$
% $C=\alpha Id_{r_1+r_2}.$

Note that, if a block $B_{12}$ has a zero-row $k$ and a zero-column $j,$ then by
 adding column $j$ to column $k$ and row $k$ to row $j$ we construct a
nonscalar endomorphism, hence $B$ is not simple. In particular, if $r_1=r_2$ then $B_{12}$
is a square matrix and can be reduced to the identity matrix $I.$
Having $B_{12}=I$
we can reduce one of matrices  $B_1$ and $B_2,$ let us say $B_1,$
to  zero and the other one
$B_2$ to its Jordan normal form.
If $r_2=1$ then $B_2= \framebox{$\lambda$},$ $\lambda\in\ko,$ in this case $B$ is simple, but
for $r_2>1$ the Jordan normal form has
an endomorphism, which can be extended to an endomorphism of $B.$
Therefore, if $B$ is simple then  $B_{12}$
can be reduced to  one of the following forms
\begin{align*}
\setlength\arrayrulewidth{0.1mm}
\setlength\doublerulesep{0.1mm}\doublerulesepcolor{black}
%\label{B01}
B_{12}=
\left\{
\begin {array}{ccc}
%\left(
%\begin {smallmatrix}
\begin {array}{@{}c@{} c @{}c@{} }
\hline
\vr&0 &\vr \\
\vr&I_{r_2}  &\vr\\
\hline
\end{array}
%\end {smallmatrix}
%\right)
&
\hbox{if}& r_1> r_2,\\
\\
\begin {array}{@{}c@{} c @{}c@{} c  @{}c@{}}
\hline
\vr~ I_{r_1} &&  \mbox{0\,} &\vr\\
\hline
\end{array} \phantom{R}
&\hbox{if}& r_2>r_1,\\
\\
\framebox{1}
&\hbox{if}& r_1=r_2=1,\\
\end{array}
\right .
\end{align*}
From the system of equations (**) we get that in case $r_1>r_2$
 block $B_2$ can be reduced to the zero matrix
and block
$B_1$ to the upper triangular block-matrix formed by three nonzero subblocks
$(B_{1.1}, B_{1.12}, B_{1.2}).$
Long but straightforward calculations show that the transformations
of $B$ which preserve already reduced blocks are uniquely determined by the
automorphisms of $B_1$ in the category  $\BM_s.$ Moreover,
$\End_{\BM_s}(B_1)=\End_{\BM_s}(B).$

In the same way the matrix $B$ can be reduced in case $r_2>r_1.$
Thus the problem $\BM_s(r_1,r_2)$ is {\it self-reproducing,} that means  we get a bijection
between
 $\BM_s(r_1,r_2)$ and $\BM_s(r_1-r_2,r_2)$ if $r_1>r_2,$
 between $\BM_s(r_1,r_2)$ and $\BM_s(r_1,r_2-r_1)$ if $r_2>r_1,$ and if $r_1=r_2>1$ then $\BM_s(r_1,r_1)$ is empty.
 In this reduction one can easily recognize the Euclidian algorithm.
Moreover,  the reduction terminates after finitely many steps when we achieve $r_1=r_2=1.$
Without loss of generality we may assume that the matrix $B\in \BM_s(1,1)$ has the form

$$
%\begin{align}
B=
\setlength\arrayrulewidth{0.1mm}
\setlength\doublerulesep{0.1mm}\doublerulesepcolor{black}
\begin {array}{@{}c@{} c @{}c@{} c @{}c@{}}
\hline\hline
\fvr& \mk{0}&\fvr&\mk{1} &\fvr\\
\hline\hline
\fvr&\mk{} &\fvr& \mk{\lambda}&\fvr\\
\hline\hline
\end{array}.
%\label{form11}
%\end{align}
\eqno{(***)}
$$
(Note that this matrix can be equivalently writen   as 
$
\begin{smallmatrix}
\hline
\vrule&^\lambda \rule{0pt}{8pt} &\vrule& ^ 1&\vrule\\
\hline
\vrule &\rule{0pt}{8pt}  &\vrule &^0&\vrule\\
\hline
\end{smallmatrix}). 
$

Objects of $\BM_s(1,1)$ are parametrized by a continuous parameter
 $\lambda \in \ko,$ thus the same holds for
$\BM_s(r_1,r_2)$ with coprime $r_1$ and $r_2.$

Let $\kE$ be a vector bundle of rank $r$ and degree $d$ with normalization
$\bn=\on(c)^{r_1}\oplus \on(c+1)^{r_2}.$
Taking into account that by Riemann-Roch theorem (Theorem \ref{RR})
$r_2=d\bmod r $ and $r_1+r_2=r ,$
we obtain  the statements about vector bundles of Theorem \ref{01}.
Moreover, if coprime integers $r>0$ and $d$
are given then for $\lambda\in \ko$ one can construct the matrix $B(\lambda)\in \BM_s(r_1,r_2),$
and hence, the unique vector bundle $\kE(r,d, \lambda),$
by reversing the reduction procedure described above:

\begin{algorithm}
Let $(r,d)\in\mZ^2$ be coprime with  positive $r$, and $\lambda\in \ko.$
\begin{itemize}
\item First,
by the Euclidean algorithm we find integers $c,$ $r_1$ and $r_2,$  $0<r_1\leq r,~~ 0\leq r_2< r$
such that $cr+r_2=d$ and $r_1+r_2=r.$ Thus we recover
the normalization sheaf $\fn=\on(c)^{r_1}\oplus \on(c+1)^{r_2}.$
\item
If $r_1=r_2=1$ the matrix $B(\lambda)$ has form (***).
\end{itemize}
Using this input data we construct the matrix $B(\lambda)\in \BM_s(r_1,r_2)$
inductively:
\begin{itemize}
\item
Let $r_1+r_2>2$ %then we proceed inductively, namely,
and $r_1>r_2.$
Assume we have the matrix \\
$B_1(\lambda)\in\BM_s(r_1-r_2, r_2),$ then $B(\lambda)\in \BM_s(r_1, r_2)$  has form
\begin{align*}
B(\lambda)=
\setlength\arrayrulewidth{0.1mm}
\setlength\doublerulesep{0.1mm}\doublerulesepcolor{black}
\begin {array}{@{}c@{}  c @{}c@{} c @{}c@{}}
\hline\hline
\fvr&  B_1(\lambda)  &\fvr&
\begin {array}{@{}c@{} c @{}c@{} }
&0 & \\
&I_{r_2}  &\\
\end{array}
&\fvr\\
\hline\hline
\fvr&\mk{} &\fvr& \mk{0}&\fvr\\
\hline\hline
\end{array}.
\end{align*}
\item
Let $r_1+r_2>2$ and $r_1<r_2.$
Assume that we have the matrix
$B_2(\lambda)\in\BM_s(r_1, r_2-r_1),$ then $B(\lambda)\in \BM_s(r_1, r_2)$  has form
\begin{align*}
B(\lambda)=
\setlength\arrayrulewidth{0.1mm}
\setlength\doublerulesep{0.1mm}\doublerulesepcolor{black}
\begin {array}{@{}c@{}  c   @{}c@{}   c    @{}c@{}}
\hline\hline
\fvr& 0      &\fvr& I_{r_1}~ 0 &\fvr\\
\hline\hline
\fvr&\mk{}   &\fvr&
\makebox[27pt][c]{\rule{0pt}{22pt}~$ ^{\displaystyle B_2(\lambda)} $ }&\fvr\\
\hline\hline
\end{array}.
\end{align*}
\item Finally, we get the matrix $\tilde{i}=i(0)+\varepsilon i_\varepsilon(0)=I_r+\varepsilon B(\lambda).$

\end{itemize}
\end{algorithm}

Let us illustrate this with a small example:
\begin{example}
Let $\kE\in \VB_s(7,12)$ be an indecomposable vector bundle of rank $7$ and degree
$12$. 
To obtain the matrices $i(0)$ and $i_\varepsilon(0)$
we calculate the normalization sheaf $\bn$ first:
$\bn=\on(1)^2\oplus \on(2)^5.$
Thus, in our notations $r_1=2$ and $r_2=5.$
The Euclidian algorithm applied to the pair $(2,5)$ gives:
%$$(2,7)\mo(2,5)\mo(2,3)\mo(2,1)\mo(1,1)$$
$$(2,5)\mo(2,3)\mo(2,1)\mo(1,1).$$
Reversing this sequence, by the above reduction procedure,
 we obtain a sequence of bijections:
$$\BM_s(1,1)\stackrel{\sim}{\longrightarrow}\BM_s(2,1)\stackrel{\sim}{\longrightarrow}
\BM_s(2,3)\stackrel{\sim}{\longrightarrow}\BM_s(2,5),$$
%$$\BM_s(1,1)\mo\BM_s(2,1)\mo\BM_s(2,3)\mo\BM_s(2,5)\mo\BM_s(2,7),$$
and finally for the matrices we get:
\vspace{-0.5cm}
$$
\setlength\arrayrulewidth{0.1mm}
\setlength\doublerulesep{0.1mm}\doublerulesepcolor{black}
\begin{array}{
@{}c@{}     c       @{}c@{}     c       @{}c@{}}
\hline\hline
%    & \hr && \hr \\
\fvr  & {\mk{0}}  &\fvr&    {\mk{1}}   & \fvr     \\
\hline\hline
\fvr&               &\fvr &        {\mk{\lambda}} & \fvr \\
%\cline{3-4}
\hline\hline
\end{array}
~\rightarrow~
\begin{array}{
@{}c@{}     c       @{}c@{}     c       @{}c@{}c       @{}c@{}}
\hline\hline
%    & \hr && \hr \\
\fvr & {\bmk{0}}  &\vr&    {\bmk{1}}   & \fvr  &{\mk{0}}   & \fvr    \\
\hline
\fvr & {\mk{0}}            &\vrule&        {\bmk{\lambda}} & \fvr  &{\mk{1}}   & \fvr\\
\hline\hline
\fvr&               &&                              & \fvr  &{\mk{0}}   & \fvr\\
%\cline{5-6}
\hline\hline
\end{array}
~\rightarrow~
\begin{array}{
@{}c@{}     c       @{}c@{}     c       @{}c@{}     c       @{}c@{}     c       @{}c@{}c       @{}c@{}}
\hline\hline
\fvr & \mk{0}  &\vr&    {\mk{0}}   &\fvr & {\mk{1}}  &\vr&    {\mk{0}}
& \vr  &{\mk{0}}   & \fvr    \\
\hline
\fvr & \mk{0}  &\vr&    {\mk{0}}   &\fvr & {\mk{0}}  &\vr&    {\mk{1}}
& \vr  &{\mk{0}}   & \fvr    \\
\hline\hline
%    & \hr && \hr \\
\fvr&&&&\fvr & {\bmk{0}}  &\vr&    {\bmk{1}}   & \vr  &{\bmk{0}}   & \fvr    \\
\cline{5-10}
\fvr&&&&\fvr & {\bmk{0}}  &\vr&    {\bmk{\lambda}} & \vr  &{\bmk{1}}   & \fvr\\
\cline{5-10}
\fvr&&&& \fvr& {\mk{0}} &\vr&{\mk{0}} &\vr   &{\bmk{0}}   & \fvr\\
%\cline{5-10}
\hline\hline
\end{array}
~\rightarrow ~
\begin{array}{
@{}c@{}     c       @{}c@{}     c  @{}c@{}     c       @{}c@{}     c       @{}c@{}     c       @{}c@{}     c       @{}c@{}c       @{}c@{}}
\hline\hline
\fvr & \mk{0}  &\vr&    {\mk{0}}   &\fvr & {\mk{1}}  &\vr&    {\mk{0}}
& \vr  &{\mk{0}}   & \vr  &{\mk{0}} & \vr  &{\mk{0}} & \fvr    \\
\hline
\fvr & \mk{0}  &\vr&    {\mk{0}}   &\fvr & {\mk{0}}  &\vr&    {\mk{1}}
& \vr  &{\mk{0}}  & \vr  &{\mk{0}} & \vr  &{\mk{0}}  & \fvr    \\
\hline\hline
\fvr&&&&\fvr & \bmk{0}  &\vr&    {\bmk{0}}   &\vr & {\bmk{1}}  &\vr&    {\bmk{0}}
& \vr  &{\mk{0}}   & \fvr    \\
\cline{5-15}
\fvr&&&&\fvr & \bmk{0}  &\vr&    {\bmk{0}}   &\vr & {\bmk{0}}  &\vr&    {\bmk{1}}
& \vr  &{\bmk{0}}   & \fvr    \\
\cline{5-15}
%    & \hr && \hr \\
\fvr&&&&\fvr&{\mk{0}}&\vr&{\mk{0}}&\vr & {\bmk{0}}  &\vr&       {\bmk{1}}   & \vr  &{\bmk{0}}
& \fvr    \\
\cline{5-15}
\fvr&&&&\fvr&{\mk{0}}&\vr&{\mk{0}}&\vr & {\bmk{0}}  &\vr&       {\bmk{\lambda}} & \vr  &{\bmk{1}}
 & \fvr\\
\cline{5-15}
\fvr&&&&\fvr&{\mk{0}}&\vr&{\mk{0}}& \vr& {\bmk{0}}  &\vr&       {\bmk{0}} &\vr   
&{\bmk{0}}
& \fvr\\
%\cline{5-10}
\hline\hline
\end{array}.
$$
\end{example}
\vspace{0.5cm}

The reduction for torsion free but not locally free sheaves can be done in a similar way.
The only difference is that the matrices $i(0)$ and $i_\varepsilon (0)$ are no longer square:
$$
i(0)=
\begin {array}{|c @{}|@{}  c  |@{}  c|}
\hline
\mk{I_{r_1}}&\mk{}& \mk{}\\
\hline
\mk{}&\mk{I_{r_2}} & \mk{}\\
\hline
\end{array}
~~\hbox{and  }~~
i_{\varepsilon}(0)=
\begin {array}{|c @{}|@{}  c  |@{}  c|}
\hline
\mk{B_1}&\mk{~B_{12}}& \mk{~B_{13}}\\
\hline
&\mk{B_{2}} &~B_{23}\\
\hline
\end{array}.
$$
The matrix $i_{\varepsilon}(0)$
has two additional blocks $B_{13}$ and $B_{23}$ with a new column size $r_3>0.$
Investigating such matrices inductively for simple matrices, we get  $r_3=1.$
Moreover, if
$r_1+1$ and $r_2+1$ are coprime then there is a unique simple matrix $\tilde{i},$
and there is no simple matrices otherwise. This unique simple matrix $\tilde{i}$
corresponds to the unique torsion free but not locally free sheaf $\ff,$
which
can be considered
as a {\it compactifying} object of the family $\VB_s(r,d).$
Let us illustrate this for small ranks:

\begin{example}
Vector bundles $\kE$ of $\VB_s(1,0)$
have $\widetilde{\kE}=\on$ as the normalization sheaf  and the corresponding
matrices $\tilde{i}$ are $\framebox{1} + \varepsilon\framebox{$\lambda$},$
 $\lambda\in \ko.$
For the unique torsion free but not locally free sheaf $\ff$
of rank $1$ and degree $0,$
 one computes that $deg(\fn)=deg(\ff)-1=-1,$
 thus $\fn=\on(-1)$ and the corresponding matrix $\tilde{i}$ is
$$
\begin{array}{
@{}c@{}     c       @{}c@{}     c       @{}c@ {}}
\hline
\vr&     {\mk{1}}    &\vr&    {\mk{0}}   &\vr     \\
\hline
\end{array}
+
\varepsilon\cdot
\begin{array}{
@{}c@{}     c       @{}c@{}     c       @{}| }
\hline
\vrule&     {\mk{0}}    &\vrule&    {\mk{1}}        \\
\hline
\end{array}.
$$
\end{example}

\begin{example}
Vector bundles $\kE$ from $\VB_s(2,1)$
have as normalization sheaf
 $\widetilde{\kE}= \on\oplus\on(1)$
thus the corresponding matrices are
$$
\tilde{i}=
\setlength\arrayrulewidth{0.1mm}
\setlength\doublerulesep{0.1mm}\doublerulesepcolor{black}
\begin{array}{
@{}c@{}     c       @{}c@{}     c       @{}c@{}  }
\hline\hline
\fvr&     {\mk{1}}    &&  {\mk{0}}    & \fvr\\
\hline\hline
\fvr&     {\mk{0}} &&  {\mk{1}}      &\fvr  \\
\hline\hline
\end{array}
\begin{array}{@{}c@{}}
{~\scriptstyle 0} \mk{}\\
{~\scriptstyle 1}\mk{}
\end{array}
+
\varepsilon\cdot
\begin{array}{
@{}c@{}     c       @{}c@{}     c       @{}c@{} }
\hline\hline
\fvr&     {\mk{0}}  &&  {\mk{1}}    &\fvr    \\
\hline\hline
\fvr&     {\mk{0}} &&  {\mk{\lambda}}       &\fvr \\
\hline\hline
\end{array}
\begin{array}{@{}c@{} c}
&{~\scriptstyle 0} \mk{}\\
&{~\scriptstyle 1}\mk{,}
\end{array}
$$
where $\lambda\in \ko.$
The normalization sheaf $\fn$ of  the torsion free but not locally free sheaf $\ff\in \TF_s(2,1)$
 has degree
$\deg(\fn)=\deg(\ff)-1=0$ thus $\fn=\on^2$ and the corresponding matrix is
$$
\tilde{i}=
\setlength\arrayrulewidth{0.1mm}
\setlength\doublerulesep{0.1mm}\doublerulesepcolor{black}
\begin{array}{
@{}c@{}     c       @{}c@{}     c       @{}c@{}     c       @{}c@{} }
\hline\hline
\fvr  &     {\mk{1}}    &&  {\mk{0}} &\fvr&        {\mk{0}}  &\fvr   \\
%\hline
\fvr &     {\mk{0}}    &&  {\mk{1}}    &\fvr&        {\mk{0}}  &\fvr   \\
\hline\hline
\end{array}
\begin{array}{@{}c@{} c}
&{~\rule{0pt}{12pt} {\scriptstyle 0}} \\
\end{array}
+
\varepsilon\cdot
\begin{array}{
@{}c@{}     c       @{}c@{}     c       @{}c@{}     c       @{}c@{} }
\hline\hline
\fvr&     {\mk{0}} &&  {\mk{1}}    &\fvr&        {\mk{0}}    &\fvr \\
%\hline
\fvr&     {\mk{0}} &&  {\mk{0}} &\fvr&        {\mk{1}}       &\fvr \\
\hline\hline
\end{array}
\begin{array}{@{}c@{} c}
&{~\rule{0pt}{12pt} {\scriptstyle 0}} \\
\end{array}
.
$$

\end{example}
\vspace{0.7cm}

\subsection{Coherent sheaves on degenerations of
elliptic curves and  Fourier-Mukai transforms}\label{subsec:FMT}

The technique of Fourier-Mukai transforms on elliptic curves led to a classification of
indecomposable coherent sheaves. This can be generalized to  the case of singular
Weierstra\ss{} cubic curves.

\begin{theorem}[\cite{BurbanKreussler1}]
Let $\EE$ be an irreducible projective curve of arithmetic genus one over an algebraically
closed field $\kk$, $p_0 \in \EE$ a
smooth point and
$\kP = \kI_\Delta \otimes \pi_1^*(\kO(p_0)) \otimes  \pi_2^*(\kO(p_0))$, where $\kI_\Delta$ is
the ideal sheaf of the diagonal $\Delta \subset \EE\times \EE$.
We have the following properties of the Fourier-Mukai transform
 $\Phi = \Phi_{\mathcal{P}}$:
\begin{enumerate}
\item $\Phi$ is an exact equivalence and $\Phi\circ\Phi \cong i^*[-1]$, where
$i: \EE\lar \EE$ is an involution of $\EE$.
\item $\Phi$ transforms semi-stable sheaves to semi-stable ones and stable sheaves to
stable ones.
\item In particular, $\Phi$ induces an equivalence between the
abelian categories\\ $\Coh^\nu(\EE)$ and $\Coh^{-\frac{1}{\nu}}(\EE)$, where
$\nu \in \mathbb{Q}\cup\{\infty\}$.
\item Let $\kF$ be a semi-stable  sheaf  of degree zero. Then the  sequence
$$
0 \lar H^0(\kE(p_0))\otimes \kO \stackrel{ev}\lar \kE(p_0) \lar coker(ev) \lar 0
$$
is exact. Moreover, the functor $\Phi(\kE) \cong  coker(ev)$ establishes
 an equivalence between the category $\Coh^0(\EE)$ of semi-stable torsion-free sheaves
of degree zero
and the category of torsion  sheaves $\Coh^\infty(\EE)$.
\end{enumerate}
\end{theorem}

From this theorem follows that for any pair of coprime integers
 $(r, d) \in \mathbb{Z}^2, r >0$
the moduli space $M_{\EE}(r, d)$ of  stable sheaves of rank $r$ and degree
$d$ is isomorphic to  $\EE$.   The unique singular point of $M_{\EE}(r, d)$
corresponds to the stable sheaf which is not locally  free.

Let
$\kT$ be an indecomposable
torsion sheaf on $\EE$.  If $\kT$ has support at a smooth point $p \in \EE$
than
$\kT \cong \kO_{\EE,p}/\mathfrak{m}_p^n$ for some $n > 0$.

The structure of torsion sheaves supported at the singular point
$s \in \EE$  is much more complicated.
First of all note  that the categories  of finite-dimensional modules over $\kO_{\EE,s}$ and
$\widehat{\kO_{\EE,s}}$ are equivalent. So, in order to understand semi-stable sheaves
on singular Weierstra\ss{} curves we have to analyze the structure of finite-dimensional
representations of $\kk\llbracket x,y\rrbracket  /(xy)$ and $\kk\llbracket x,y\rrbracket  /(y^2 - x^3)$ first.

Let $\mathbf{R} = \kk\llbracket x,y\rrbracket  /(xy)$, then it is easy to show that
all indecomposable
finite length $\mathbf{R}$-modules  generated by one element are
$\mathcal{M}((n,m), 1,\lambda) = \mathbf{R}/(x^n + \lambda y^m)$ for  $n,m \ge 1$,
$\lambda \in \kk^*$ and  $\mathcal{N}(0,(n, m),0) = \mathbf{R}/(x^{n+1}, y^{m+1})$ for
 $n,m \ge 0$.
A classification of all indecomposable $\mathbf{R}$--modules  was obtained  by
Gelfand and Ponomarev
\cite{GP} and independently by Nazarova and Roiter \cite{NazRoi},
see also \cite{BurbanDrozd1} for a description via derived categories.
We identify an indecomposable  torsion module $\kT$ supported at $s$ with the corresponding
$\kk\llbracket x,y\rrbracket  /(xy)$--module.

\begin{theorem}[\cite{BurbanKreussler1}]
The Fourier-Mukai  transform  $\Phi_\kP$ maps the torsion module\\
$\mathcal{M}((n,m), 1,\lambda)$ to the degree zero
 semi-stable vector bundle
$$\kB((\overbrace{1,0,\ldots, 0}^{m},
\overbrace{-1, 0,\ldots,0}^{n}), 1, (-1)^{(n+m)}\lambda)$$   and
$\mathcal{N}(0,(n, m),0)$ to the semi-stable
torsion  free sheaf $$\kS(\overbrace{0,\ldots, 0}^m,-1,\overbrace{0,\ldots, 0}^n).$$
\end{theorem}

In \cite{BurbanKreussler1} a complete correspondence between torsion sheaves and
semi-stable sheaves of degree zero was described. Using the
 technique of  relative Fourier-Mukai transforms one gets   a powerful tool to construct
interesting examples of relatively stable and semi-stable sheaves on elliptically fibered
varieties, see  \cite{FMW, BurbanKreussler2}.

In a similar way, if $\mathbf{R} = \mathbf{k}\llbracket x,y\rrbracket  /(y^2 - x^3)= \kk\llbracket t^2, t^3\rrbracket  $,
then a finite length $\mathbf{R}$--module is given by a finite
dimensional vector space $V$ over $\kk$ and two endomorphisms $X,Y : V\rightarrow V$ which
satisfy $Y^2 - X^3 = 0$.
 It is again very easy to
classify all $\mathbf{R}$--modules of  the form  $\mathbf{R}/I$,
where $I$ is an ideal in $\mathbf{R}$:
there are  one-parameter families  of modules of projective dimension one:
$\mathbf{R}/(t^n + \lambda t^{n+1})$,  $n \ge 2$  
 and  $ \lambda \in \kk$,  and
discrete  series of modules of infinite projective dimension,
$\mathbf{R}/(t^n,t^{n+1})$, where  $n\ge 2$.

However,
there is an essential difference between the rings $\kk\llbracket x,y\rrbracket  /(xy)$ and
$\mathbf{R} = \kk\llbracket x,y\rrbracket  /(y^2 - x^3)$: the first  has tame representation
type \cite{GP, NazRoi} whereas the second  is  wild.

\begin{proposition}[\cite{Drozd}]
The category of finite length modules over the ring \\ $\kk\llbracket x,y\rrbracket  /(y^2-x^3)$ is wild.
\end{proposition}

\noindent
\emph{Proof}.
We have a fully-faithful exact functor
$
\kk\langle z_1,z_2\rangle \lar \mathsf{Rep}_\Gamma,
$
where $\Gamma$ is the quiver
$\xymatrix{\ar@(ul,dl)_{a} \bullet \ar[r]^b & \bullet}$. This functor  maps
a $\kk\langle z_1,z_2\rangle$--module $(Z_1, Z_2, \kk^n)$ to the representation of $\Gamma$ given by
$$
A =
\left(
\begin{array}{cc}
Z_1 & Z_2 \\
I_n   & 0
\end{array}
\right),\quad
B = \left(
\begin{array}{cc}
0 & I_n \\
\end{array}
\right).
$$
Moreover,
we have another fully-faithful exact functor
$$\mathsf{Rep}_\Gamma \lar \kk\llbracket x,y\rrbracket  /(y^2 - x^3),$$
mapping a representation  $(A,B, \kk^n, \kk^m)$ to the $\kk\llbracket x,y\rrbracket  /(y^2 - x^3)$--module given by the
matrices
$$
Y = \left(
\begin{array}{ccc}
0_{3n} & 0 & I \\
0 & 0_m & 0 \\
0 & 0 & 0_{3n}
\end{array}
\right),\quad
X = \left(
\begin{array}{ccc}
X_1 & 0 & X_2 \\
0 & 0 & X_3 \\
0 & 0 & X_1
\end{array}
\right),
$$
where
$$
X_1 = \left(
\begin{array}{ccc}
0_n & 0 & I \\
0 & 0_n & 0 \\
0 & 0 & 0_n
\end{array}
\right), \quad
X_2 = \left(
\begin{array}{ccc}
0_n & 0 & 0   \\
I & 0_n & 0 \\
0 & A & 0_n
\end{array}
\right),  \quad
X_3 =
%\left(
%\begin{array}{ccc}
(0_{m\times n} ~ B_{m\times n} ~ 0_{m\times n}).
%\end{array}
%\right).
$$
Taking the composition of these two functors we see that the category of finite dimensional
$\kk\llbracket x,y\rrbracket  /(y^2-x^3)$--modules is wild.

\medskip

Since via an appropriate Fourier-Mukai transform the category of torsion modules
over the ring $\widehat{\kO_{\EE,s}}$
is equivalent to the category of semi-stable
torsion free sheaves of a given slope $\nu$ with  non-locally free Jordan-H\"older quotient, we
obtain  the following corollary.

\begin{corollary}\label{cor:cuspiswild}
Let $\mathbb{E}$ be a cuspidal Weierstra\ss{} curve and
$\nu\in \mathbb{Q}\cup\{\infty\}$, then
the category $\Coh^\nu(\mathbb{E})$ of semi-stable torsion free sheaves of slope $\nu$
on $\mathbb{E}$
has wild  representation type.
\end{corollary}

%%%%%%%%%%%%%%%%%%%%%%%%%%%%%%%%%%%%%%%%%%%%%%%%%%%%%%%%%%%%%%%%%%%%%%

%%%%%%%%%%%%%%%%%%%%%%%%%%%%%%%%%%%%%%%%%%%%%%%%%%%%%%%%%%%%%%%%%%%%%%

\printindex
\end{document}